\documentclass[11pt,a4]{article}
\usepackage[T2A]{fontenc}
\usepackage[cp866]{inputenc}
\usepackage[mathscr]{eucal}
\usepackage{amssymb}
\usepackage{amsmath}
\renewcommand{\leq}{\leqslant}
\renewcommand{\geq}{\geqslant}

\newcommand{\R}{{\mathbb R}}

\renewcommand{\k}{\rule{0.7em}{0.7em}}
\usepackage{fancyhdr}
\begin{document}
\sloppy

{\normalsize

\thispagestyle{empty}

\mbox{}
\\[-2.25ex]
\centerline{
{
\bf
COMPACT INTEGRAL MANIFOLDS OF DIFFERENTIAL SYSTEMS
}
}
\\[2.75ex]
\centerline{
\bf
V.N. Gorbuzov
}
\\[1.75ex]
\centerline{
\it
Department of Mathematics and Computer Science,
Yanka Kupala Grodno State University,
}
\\[1ex]
\centerline{
\it
Ozeshko 22, Grodno, 230023, Belarus
}
\\[1ex]
\centerline{
E-mail: gorbuzov@grsu.by
}
\\[5.75ex]
\centerline{{\large\bf Abstract}}
\\[1ex]
\indent
The boundedness tests for the number of compact integral manifolds of autonomous
ordinary differential systems, of  autonomous  total differential systems, of linear systems of partial
differential equations, of Pfaff systems of equations, and of
systems of exterior differential equations are proved.
\\[2ex]
\\[1.25ex]
\\[5.5ex]
\centerline{{\large\bf Contents}}
\\[1.5ex]
{\bf  Introduction}\footnote[1]{
The basic results of this paper have been published in the monograph
Gorbuzov V.N. {\it Integrals of  differential systems} (Russian), Grodno State University, Grodno, 2006.
}
                                                                                                                   \hfill\ 2
\\[0.5ex]
{\bf
1. The boundedness of the number of compact regular integral manifolds}\footnote[2]{
The initiative version of this Section has been published in
{\it Differentsial Uravneniya} ({\it Differential Equations}), Vol. 35 (1999), No. 10, 1325-1329;
Vol. 36 (2000), No. 11, 1563-1565.
}
                                                                                                                    \hfill\ 4
\\[0.5ex]
\mbox{}\hspace{1em}
1.1. Autonomous ordinary differential system                                                                                             \dotfill\ 4
\\[0.5ex]
\mbox{}\hspace{1em}
1.2. Autonomous total differential system                                                     \dotfill\ 8
\\[0.5ex]
\mbox{}\hspace{2.75em}
1.2.1.\!\! Test of the boundedness of the number of compact regular
integral manifolds
                                                                                                                         \dotfill\ 8
\\[0.5ex]
\mbox{}\hspace{2.75em}
1.2.2.\!
Tests of the absence of compact regular orbits
                                                                                                                         \dotfill\ 9
\\[1ex]
\noindent
{\bf
2. The boundedness of the number of compact integral hypersurfaces}\footnote[3]{
The results of this Section were originally published in
{\it Differential Uravneniya {\rm(}Differential Equations}), Vol. 33 (1997), No. 10, 1307-1311;
Vol. 35 (1999), No. 1, 30-37.
}
                                                                                                                              \hfill\ 11
\\[0.5ex]
\mbox{}\hspace{1em}
2.1. System of exterior differential equations
                                                                                                                                   \dotfill\ 11
\\[0.5ex]
\mbox{}\hspace{1em}
2.2. Pfaff system of equations
                                                                                                                                   \dotfill\ 14
\\[0.5ex]
\mbox{}\hspace{2.75em}
2.2.1.\! The first boundedness test for the number of compact integral hypersurfaces
                                                                                                                                    \dotfill\ 14
\\[0.5ex]
\mbox{}\hspace{2.75em}
2.2.2.\! The second boundedness test for the number of compact
\\
\mbox{}\hspace{5.5em}
integral hypersurfaces
                                                                                                                                    \dotfill\ 16
\\[0.5ex]
\mbox{}\hspace{1em}
2.3. Ordinary autonomous differential system
                                                                                                                                   \dotfill\ 19
\\[0.5ex]
\mbox{}\hspace{2.75em}
2.3.1.\! The first boundedness test for the number of compact integral hypersurfaces
                                                                                                                                    \dotfill\ 19
\\[0.5ex]
\mbox{}\hspace{2.75em}
2.3.2.\! The second boundedness test for the number of compact
\\
\mbox{}\hspace{5.5em}
integral hypersurfaces
                                                                                                                                    \dotfill\ 20
\\[0.5ex]
\mbox{}\hspace{2.75em}
2.3.3.\! Test of the boundedness of the number of compact regular
\\
\mbox{}\hspace{5.5em}
integral hypersurfaces
                                                                                                                                    \dotfill\ 22
\\[0.5ex]
\mbox{}\hspace{1em}
2.4. Autonomous total differential system
                                                                                                                                   \dotfill\ 23
\\[0.5ex]
\mbox{}\hspace{1em}
2.5. Linear homogeneous system of partial differential equations
                                                                                                                                   \dotfill\ 25
\\[0.5ex]
{\bf  References}                                                                                                  \dotfill\ 26

\newpage

\mbox{}
\\
\centerline{
\large\bf
Introduction
}
\\[1.5ex]
\indent
At the Second International Congress of Mathematicians, which met at Paris in 1900,
the outstanding German mathematician David Hilbert made the report "Mathematical Problems"
[1, pp. 290 -- 329]. In it, he presented 23 problems, extending over all fields of mathematics.
Hilbert's sixteenth problem is the investigation of the topology of algebraic curves and surfaces.
The second part of this problem is the problem on the maximal number and the arrangement of isolated closed trajectories
(limit cycles [2]) of an ordinary autonomous second-order differential system with polynomial right-hand sides.
\vspace{0.25ex}

Poincar\'{e} proved  [3, pp. 112 -- 136] a test of the absence of closed trajectories (periodic solutions)
for a differential system
\\[1.75ex]
\mbox{}\hfill
$
\dfrac{dx}{dt}=P(x,y),
\qquad
\dfrac{dy}{dt}=Q(x,y)
$
\hfill (0.1)
\\[2ex]
where $\!P\colon G\!\to\! \R\!$ and $\!Q\colon G\!\to\! \R\!$ are
\vspace{0.75ex}
continuously differentiable functions on a domain $\!G\!\subset\! \R^2\!.$

{\bf The Poincar\'{e} test.}
{\it
If there exists the continuously differentiable on a domain $G^{\;\!\prime}\subset G$
function $N_2^{}\colon G^{\;\!\prime}\to \R$ such that the function
\\[1.5ex]
\mbox{}\hfill
$
H_2\colon (x,y)\to\, \partial_x N_2(x,y)\;\!P(x,y)+ \partial_y N_2(x,y)\;\!Q(x,y)
$
\ for all
$
(x,y)\in G^{\;\!\prime}
\hfill
$
\\[1.5ex]
is definite, then the system {\rm (0.1)} has no closed trajectories in the domain $G^{\;\!\prime}.$
}
\vspace{0.5ex}

Bendixson proved [4] a test for the absence of simple closed curves, which
are made from trajectories of system (0.1). This was the first test based on the divergence of the vector field of
the differential system (0.1) and proved with the help of the Green formula
about relations of curvilinear and double integrals.

A simple closed curve, which is made from trajectories of system (0.1), either is a closed trajectory or
is made from alternation among themselves of whole nonclosed trajectories  and equilibria.
In particular a limit cycle is such curve.
\vspace{0.75ex}

{\bf The Bendixson test}.
{\it
If the divergence of the vector field
\\[1.5ex]
\mbox{}\hfill
$
a\colon (x,y)\to (P(x,y),Q(x,y))
$
\ for all $(x,y)\in G
\hfill
$
\\[1.5ex]
on a simply connected domain $G^{\;\!\prime}\subset G$ is constant sign, then
the system {\rm (0.1)} in the domain $G^{\;\!\prime}$ has no simple closed curves, which
are made from trajectories of system {\rm (0.1)}.
}
\vspace{0.75ex}

Bendixson's test was generalized by Dulac [5] and generalization is referred to as the Bendixson --- Dulac test
for the absence in a simply connected domain
of simple closed curves, which are made from trajectories of system {\rm (0.1)},
and for the existence of at most one simple closed curve, which is made from trajectories of system {\rm (0.1)},
in an annular domain.
\vspace{0.75ex}

{\bf The Bendixson --- Dulac test.}
{\it
If there exists the continuously differentiable on a simply connected domain $G^{\;\!\prime}\subset G$
function $\mu \colon G^{\;\!\prime}\to \R$  such that the divergence of the vector field
\\[1.5ex]
\mbox{}\hfill
$
b\colon (x,y)\to\, \mu(x,y)\;\!(P(x,y),Q(x,y))$
\ for all
$
(x,y)\in G^{\;\!\prime}
\hfill
$
\\[1.5ex]
is  constant sign, then
the system {\rm (0.1)} in the domain $G^{\;\!\prime}$ has no simple closed curves, which
are made from trajectories of system {\rm (0.1)}.
}
\vspace{0.75ex}

{\bf The Dulac test.}
{\it
If there exists the continuously differentiable on a biconnected domain $G^{\;\!\prime}\subset G$
scalar function $\mu \colon G^{\;\!\prime}\to \R$  such that the divergence of the vector field
\\[1.5ex]
\mbox{}\hfill
$
b\colon (x,y)\to\, \mu(x,y)\;\!(P(x,y),Q(x,y))$
\ for all
$
(x,y)\in G^{\;\!\prime}
\hfill
$
\\[1.5ex]
is  constant sign, then
the ordinary differential system {\rm (0.1)} in the domain $G^{\;\!\prime}$ has
at most one simple closed curve, which is made from trajectories of system {\rm (0.1)},
such that
this curve contains the inner boundary of the domain $G^{\;\!\prime}.$
}
\vspace{0.75ex}

In [6], Bendixson's and Dulac's tests for system (0.1) was generalized on a multiply connected domain by the present author.

\newpage

{\bf Theorem 0.1}
{\sl {\rm(}a test for the boundedness of the number of limit cycles of
an ordinary autonomous second-order differential system}).
{\it
If there exists the constant sign continuously differentiable on an $s\!$-connected domain
$G^{\;\!\prime}\subset G$
scalar function $\mu \colon G^{\;\!\prime}\to \R$  such that the divergence of the vector field
$b\colon (x,y)\to \mu(x,y)(P(x,y),Q(x,y))$ for all $(x,y)\in G^{\prime}$ is
constant sign or is identically zero on the domain $G^{\prime},$
then the differential system {\rm (0.1)} in the domain $G^{\prime}$ has not more than $s-1$ limit cycles.
}
\vspace{0.35ex}

The aim of this paper is
the investigation of the topology of compact integral manifolds of differential systems:

an ordinary autonomous differential system of $n\!$-th order
\\[1.75ex]
\mbox{}\hfill             
$
\dfrac{dx}{dt}=f(x),
$
\hfill (D)
\\[2.25ex]
where
\vspace{1ex}
$t\in\R,\ x\in\R^n,$
a vector column $\dfrac{dx}{dt}={\rm colon}\, \Bigl(\dfrac{dx_1}{dt}\,,\ldots,\dfrac{dx_n}{dt}\Bigr),$
a vector function $ f(x)={\rm colon}(f_1^{}(x),\ldots,f_n^{}(x)),\!$
\vspace{0.75ex}
scalar functions $f_i^{}\colon G\to\R,\ i\!=\!1,\ldots,n,$ a domain $G\subset\R^{n};$

an autonomous system of total differential equations
\\[1.75ex]
\mbox{}\hfill             
$
dx=X(x)\,dt,
$
\hfill (TD)
\\[2ex]
where
\vspace{0.75ex}
$\!t\in\R^m,\, x\in\R^n,\, m<n,\ dx={\rm colon}\, (dx_1,\ldots,dx_n),\,
dt={\rm colon}\, (dt_1,\ldots,dt_m),\!$
an
$\!n\!\times\! m$
matrix
\vspace{1ex}
$X(x)=\bigl|\!\bigl| X_{ij}^{}(x) \bigl|\!\bigl|,\
X_{ij}^{} \colon G\to \R,\  i=1,\ldots,n,\, j=1,\ldots,m,$
a domain $G\subset\R^{n};$

a linear homogeneous system of partial differential equations
\\[1.75ex]
\mbox{}\hfill                              
$
{\frak X}_{{}_{\scriptstyle j}}(x)\;\!y = 0,
\quad
j = 1,\ldots,m,
$
\hfill  (\!$\partial$\!)
\\[1.5ex]
where $\!x\in\R^n,\, y\in\R,\, m<n,\!$
\vspace{0.35ex}
linear differential ope\-ra\-tors
of first order
$
\!{\frak X}_{{}_{\scriptstyle j}}(x)\! =\!
\sum\limits_{i=1}^{n}\!
X_{ij}^{}(x)\;\!\partial_{{}_{\scriptstyle x_i}}$ for all $x \in G,\
X_{ij}^{}\colon G\to \R,\, i=1,\ldots,n, \ j=1,\ldots,m,$
a domain $G\subset\R^{n};$
\vspace{0.75ex}

a Pfaff system of equations
\\[1.5ex]
\mbox{}\hfill                              
$
\omega_{j}^{}(x)=0,
\quad
j = 1,\ldots,m,
$
\hfill {\rm (Pf)}
\\[1.75ex]
where $x\in\R^n,\ m<n,$
\vspace{0.35ex}
linear differential forms
$
\omega_{j}^{}(x)=
\sum\limits_{i=1}^{n}
w_{ji}^{}(x)\;\!dx_i^{}$ for all
$x \in G,$
$
w_{ji}^{}\colon G\to \R,\ j=1,\ldots,m,\
i=1,\ldots,n,$
a domain $G\subset\R^{n};$
\vspace{0.5ex}

a system of exterior differential equations
\\[1.75ex]
\mbox{}\hfill                                       
$
\zeta_j^{}(x) =  0,
\quad
j = 1,\ldots,m,
$
\hfill {\rm (ED)}
\\[2ex]
where $x\in\R^n,\ m<n,$
coefficients of  $p_{{}_{\scriptstyle j}}\!$-forms
$\zeta_{j}^{}\;\!, \ 1 \leq p_{{}_{\scriptstyle j}} \leq n - 1,\ j = 1,\ldots,m,$
are scalar functions on a domain $G\subset\R^{n}.$

Let the linear differential ope\-ra\-tors
\vspace{0.25ex}
${\frak X}_{j}^{},\ j=1,\ldots,m,$ of the system of partial differential equations (\!$\partial$\!),
the linear differential forms $\omega_{j}^{},\,j=1,\ldots,m,$
\vspace{0.25ex}
of the Pfaff system of equations (Pf),
and the $p_{{}_{\scriptstyle j}}\!$-forms $\zeta_{j}^{},\ j = 1,\ldots,m,$ of the system of exterior differential equations {\rm (ED)}
be not linearly bound on the domain $G$ [7, pp. 105 -- 115].
\vspace{0.25ex}
Note also that the operators ${\frak X}_{j}^{},\,j=1,\ldots,m$
(the 1-forms $\omega_{j}^{},\ j=1,\ldots,m,$
the $p_{{}_{\scriptstyle j}}\!$-forms $\zeta_{j}^{},\ j = 1,\ldots,m)$
\vspace{0.25ex}
are called {\it linearly bound} on the domain $G$ if
these operators (1-forms, $p_{{}_{\scriptstyle j}}\!$-forms) are
linearly dependent in any point of the domain $G.$

Note that the Poincar\'{e} test was generalized by V.F. Tkachev in [8] for system {\rm (D)}.
\vspace{0.5ex}

{\bf Theorem 0.2.}
{\it
If there exists the continuously differentiable on a domain $G^{\;\!\prime}\subset G$
function $N\colon G^{\;\!\prime}\to \R$ such that the function
$
H\colon x\to
\sum\limits_{i=1}^{n}\partial_{x_i}^{}N(x)\;\!f_{i}^{}(x)
$
for all $x\in G^{\;\!\prime}$
is definite,
then the system {\rm (D)} with the
\vspace{0.35ex}
continuously differentiable on the domain $G^{\;\!\prime}$ vector function
$f\colon G\to \R^{n}$
has no closed trajectories in the domain $G^{\;\!\prime}.$
}
\vspace{0.5ex}

The basic results about the topology of compact integral manifolds of differential systems were
originally published in the papers [9, 10, 11] and then in the monograph [12].

To avoid ambiguity, we stipulate the following notions.

We recall that by domain we mean open arcwise connected set.
\vspace{0.25ex}

By $f\in C^{\;\!k}_{}(G)\ \bigl(f\in C^{\;\!\infty}_{}(G)\bigr)$ we denote
\vspace{0.25ex}
a $k\!$-times continuously differentiable (holomorphic) function on the domain $G.$
If the entries of the functional matrix $X,$
the coordinate functions of the linear differential ope\-ra\-tors ${\frak X}^{}_{j}\;\!,$ and
\vspace{0.35ex}
the coefficient of the differential forms $\omega^{}_{j}\;\!,\ \zeta^{}_{j}$ are
\vspace{0.35ex}
$k\!$-times continuously differentiable (holomorphic) on the domain $G,$ then we write
$\!\!
X\!\in\! C^{\;\!k}(G),\
{\frak X}^{}_{j}\!\in\! C^{\;\!k}(G),\
\omega^{}_{j}\!\in\! C^{\;\!k}(G),\
\zeta^{}_{j}\!\in\!  C^{\;\!k}(G)
\vspace{0.35ex}
\
\bigl(
X\!\in C^{\;\!\infty}(G),\,
{\frak X}^{}_{j}\!\in C^{\;\!\infty}(G),\,
\omega^{}_{j}\!\in C^{\;\!\infty}(G),
$
$
\zeta^{}_{j} \in C^{\;\!\infty}(G)\bigr),
\ j=1,\ldots,m.
$
\vspace{0.35ex}

If the scalar function $g\colon G\to \R$ such that either $g(x)\geq 0$ for all $x\in G$
or $g(x)\leq 0$ for all $x\in G,$ and if the relation $g(x)=0$
is possible only on a null-set with respect to the $\nu\!$-di\-men\-si\-o\-nal measure, then this function is said to be
{\it $\nu\!$-constant sign} on the domain $G\subset \R^n,$ $n>\nu,$
and we single out the cases {\it $\nu\!$-positive sign} and {\it $\nu\!$-negative sign} functions on the domain $G.$
If $\nu=n,$ then we say that the scalar function $g\colon G\to \R$ is {\it constant sign},
{\it sign-positive}, and {\it sign-negative} on the domain $G,$ respectively.
\vspace{0.35ex}

If either the function $g$ is positive $(g(x)>0$ for all $x\in G)$ or
\vspace{0.25ex}
the function $g$ is negative $(g(x)<0$ for all $x\in G),$ then we say that
the function $g$ on the domain $G$ is {\it definite}.
\vspace{0.35ex}

By a {\it gap} $\theta$ in the domain $G$ of the space $\R^n$
\vspace{0.25ex}
with homotopy group $\pi^{}_{\nu}(G),\ \nu\leq n-1,$
we mean a nonempty arcwise connected set $\theta$ such that $\theta\cap G=\O$
and in the domain $G$ there exists a manifold homeomorphic to a sphere $S^{\nu}$
such that $\theta$ is an obstruction to the continuous contraction of this manifold into a point.

An integral manifold of the system of total differential equations {\rm (TD)}\, (of the ordinary differential system (D))
is said to be {\it regular} if it is oriented and does not contain singular points (equilibria) of this differential system.

A compact integral manifold of dimension $\nu$ for the total differential system {\rm (TD)}\,
(for the ordinary differential system (D))
is called {\it isolated} if in some $\varepsilon\!$-neighborhood of this manifold there is
no other compact integral manifolds of dimension $\nu$ for this system.

In the same way, we may define isolated closed trajectories (limit cycles) of system {\rm (D)},
isolated compact regular orbits [13, 14] of the completely solvable system {\rm (TD)},
isolated compact integral hypersurfaces of systems {\rm(D)}, {\rm(TD)}, $(\partial\;\!),$ {\rm(Pf)}, and {\rm(ED)}.
\\[5.75ex]
\centerline{
\large\bf
1.\!\! The boundedness of the number of compact regular integral manifolds
}
\\[2.25ex]
\centerline{
{\bf  1.1. Autonomous ordinary differential system}
}
\\[1.75ex]
\indent
Let $\Omega_{\;\!t_{0}^{}}^{}$ be a $\nu\!$-dimensional manifold, bounded by $(\nu-1)\!$-dimensional manifolds
$\Lambda_{{}_{\scriptstyle t_0^{}}}^k,
\linebreak
k\!=\!1,\ldots,s,$ from the space $\R^n,\, 3\leq \nu \leq n.$
\vspace{0.35ex}
In addition, the manifolds $\Omega_{\;\!t_{0}^{}}^{},\, \Lambda_{{}_{\scriptstyle t_0^{}}}^k,\, k\!=\!1,\ldots,s,$
consist of points $x_0^{}=(x_1^0,\ldots,x_n^0)$ provided by solutions
\vspace{0.5ex}
$x\colon t\to x(t)$ for all $t\in J_0^{}$ of system (D) with  $f\in C^1(G)$ such that
$
x(t_0^{}) = x_0^{}, \  t_0^{}\in J_0^{},\  x_0^{}\in G,
$
\vspace{0.65ex}
and $n \geq 3.$
We define manifolds
$\Omega_{{}_{\scriptstyle t_0^{}+h}}$ and
\vspace{0.35ex}
$\Lambda_{{}_{\scriptstyle t_0^{}+h}}^k,
\, k = 1,\ldots, s,$ as sets of points $x^h = (x_1^h,\ldots,x_n^h),\,x^h\in G,$
provided by the same solutions
$x\colon t\to x(t)$ for all $t\in J_0^{}$ under
\vspace{0.75ex}
$t = t_0^{} + h,$ i.e.,
$
x^h = x(t_0^{} + h),\ x^h \in G.
$

The\-re\-fo\-re we can form the mappings
\vspace{0.5ex}
$
\Omega\colon t \to \Omega(t)$
for all $t\in {\rm U}(t_0^{})$
and $\Lambda^k\colon t \to \Lambda^k(t)$
for all $t\in {\rm U}(t_0^{}),\  k = 1,\ldots, s,$
\vspace{0.35ex}
such that
$
\Omega(t_0^{})  = \Omega_{{}_{\scriptstyle t_0^{}}}$
and
$\Lambda^k(t_0^{}) =
\Lambda^k_{{}_{\scriptstyle t_0^{}}},
\ k = 1,\ldots, s,
$
where ${\rm U}(t_0^{})$ is a neig\-h\-bor\-hood of the point $t_0^{}.$
\vspace{0.5ex}

By $\R^\nu_{\xi},$ where $\xi$ is the multiindex $\xi_1^{}\ldots\,\xi_n^{},$
we denote the subspace of the arithmetic space $\R^n$ formed
by the base coordinates $x_{\xi_1^{}}^{},\ldots,x_{\xi_\nu^{}}^{}$
\vspace{1ex}
and by $\Omega^{\;\!\xi}_{\;\!t_0^{}}$ we denote the natural projection of the manifold $\Omega_{\;\!t_0^{}}^{}$
onto the subspace $\R^{\nu}_{\xi}.$
\vspace{0.5ex}
Note that, among all $\nu\!$-dimensional subspaces $\R^{\nu}$ formed by $\nu$ coordinates of the basis
\vspace{0.5ex}
$x_1^{},\ldots, x_n^{}$ there is at least one subspace in which the manifold $\Omega^{\;\!\xi}_{\;\!t_0}$ has the dimension $\nu.$
\vspace{0.5ex}
To be definite, we assume that this is the subspace $\R^{\nu}_{\xi}.$

Let $V^{\;\!\xi}_{\;\!t_0^{}}$ be  the $\nu\!$-dimensional volume of the manifold $\Omega^{\;\!\xi}_{\;\!t_0^{}}.$
\vspace{0.5ex}
We introduce the mappings
$
\Omega^{\xi}\colon t\to \Omega^{\xi}(t)$
for all $t \in {\rm U}(t_0^{})$
and
\vspace{0.5ex}
$
V^{\xi}\colon t \to V^{\xi}(t)$
for all $t \in {\rm U}(t_0^{})$
such that
$
\Omega^{\;\!\xi}(t_0^{}) =\Omega^{\;\!\xi}_{\;\!t_0^{}},$
$V^{\;\!\xi}(t_0) =V^{\;\!\xi}_{\;\!t_0^{}}.$
In this case, the scalar function of scalar argument
\\[1.5ex]
\mbox{}\hfill
$
\displaystyle
V^{\xi}\colon t\to
\int\limits_{V^{\xi}_t}
d x^{\nu}_{\xi}
$
\ for all \
$
t \in {\rm U}(t_0),
\quad
x^{\nu}_{\xi} \in
\R^{\nu}_{\xi},
\hfill
$
\\[1.5ex]
is an additive function of the $\nu\!$-dimensional volume.
\vspace{0.5ex}

For the scalar function $V^{\xi}$ we evaluate the derivative
\vspace{0.35ex}
$
{\sf D}V^{\xi}\colon {\rm U}(t_0) \to \R,
$
specifying the change of the $\nu\!$-dimensional volume $V^{\xi}(t)$ in the neighborhood ${\rm U}(t_0^{}).$
\vspace{0.35ex}
We add the increment $\Delta t$ to the independent variable $t$ and evaluate the coefficient at $\Delta t$
\vspace{0.35ex}
in the Taylor formula for the scalar function
$
V^{\xi}\colon t + \Delta t \to V^{\xi}(t + \Delta t)$
for all $t + \Delta t \in {\rm U}(t_0^{})$
\vspace{0.35ex}
with the remainder term $o(\Delta t).$
Let $\widetilde{x}_i^{} = x_i^{}(t + \Delta t), \  i = 1,\ldots,n.$ Then
\\[2.5ex]
\mbox{}\hfill
$
\displaystyle
V^{\xi}(t+\Delta t)
=
\int\limits_{ V^{\xi}_{{}_{\scriptsize t+\Delta t}}}
d\;\!\widetilde{x}^{\,\nu}_{\xi}$
\
for all $t + \Delta t \in {\rm U}(t_0^{}),
\ \ \
\widetilde{x}^{\,\nu}_{\xi}\in \R^{\nu}_{\xi}.
\hfill
$
\\[2.5ex]
\indent
On the other hand,
\\[2.5ex]
\mbox{}\hfill
$
\displaystyle
V^{\xi}(t+\Delta t)
=
\int\limits_{V^{\xi}_t}
J\bigl( \widetilde{x}^{\,\nu}_{\xi};\, x^{\nu}_{\xi}\bigr)
\;\!dx^{\nu}_{\xi}$
\
for all $t + \Delta t \in {\rm U}(t_0^{}),
\hfill
$
\\[2.5ex]
where $J\bigl( \widetilde{x}^{\,\nu}_{\xi};\, x^{\nu}_{\xi}\bigr)$ is the Jacobian. Moreover,
\\[2ex]
\mbox{}\hfill
$
\widetilde{x}_{i}^{} =x_{i}^{} + f_{i}^{}(x)\Delta t + o(\Delta t),
\quad
 i=1,\ldots,n,
\hfill
$
\\[2ex]
and the Jacobian
\\[1.5ex]
\mbox{}\hfill
$
J\bigl( \widetilde{x}^{\,\nu}_{\xi};\, x^{\nu}_{\xi}\bigr)
 =
\hfill
$
\\[3.75ex]
\mbox{}\hfill
$
=
\left|\!\!
\begin{array}{cccc}
1+\partial_{{}_{\scriptstyle
x_{{}_{\scriptsize \xi_{{}_{\scriptsize 1}}}}}}\!
f_{{}_{\scriptstyle \xi_{{}_{\scriptsize 1}}}}\!\!(x)\Delta t + o(\Delta t)
&
\partial_{{}_{\scriptstyle x_{{}_{\scriptsize
\xi_{{}_{\scriptsize 2}}}}}}\!
f_{{}_{\scriptstyle \xi_{{}_{\scriptsize 1}}}}\!(x)\Delta t + o(\Delta t)
&
\!\!\ldots
&
\partial_{{}_{\scriptstyle x_{{}_{\scriptsize
\xi_{{}_{\scriptsize \nu}}}}}}\!
f_{{}_{\scriptstyle \xi_{{}_{\scriptsize 1}}}}\!(x)\Delta t + o(\Delta t)
\\[3.5ex]
\partial_{{}_{\scriptstyle x_{{}_{\scriptsize
\xi_{{}_{\scriptsize 1}}}}}}\!
f_{{}_{\scriptstyle \xi_{{}_{\scriptsize 2}}}}\!(x)\Delta t + o(\Delta t)
&
1+\partial_{{}_{\scriptstyle x_{{}_{\scriptsize
\xi_{{}_{\scriptsize 2}}}}}}\!
f_{{}_{\scriptstyle \xi_{{}_{\scriptsize 2}}}}\!\!(x)\Delta t + o(\Delta t)\!\!
&
\!\!\ldots
&
\partial_{{}_{\scriptstyle x_{{}_{\scriptsize
\xi_{{}_{\scriptsize \nu}}}}}}\!
f_{{}_{\scriptstyle \xi_{{}_{\scriptsize 2}}}}\!(x)\Delta t + o(\Delta t)
\\[2ex]
\cdot\ \ \ \cdot\ \ \ \cdot\ \ \ \cdot\ \ \ \cdot\ \ \ \cdot\ \ \ \cdot\ \ \ \cdot
&
\cdot\ \ \ \cdot\ \ \ \cdot\ \ \ \cdot\ \ \ \cdot\ \ \ \cdot\ \ \ \cdot\ \ \ \cdot
&
\!\!\cdot
&
\cdot\ \ \ \cdot\ \ \ \cdot\ \ \ \cdot\ \ \ \cdot\ \ \ \cdot\ \ \ \cdot\ \ \ \cdot
\\[1.25ex]
\partial_{{}_{\scriptstyle x_{{}_{\scriptsize
\xi_{{}_{\scriptsize 1}}}}}}\!
f_{{}_{\scriptstyle \xi_{{}_{\scriptsize \nu}}}}\!(x)\Delta t + o(\Delta t)
&
\partial_{{}_{\scriptstyle x_{{}_{\scriptsize
\xi_{{}_{\scriptsize 2}}}}}}\!
f_{{}_{\scriptstyle \xi_{{}_{\scriptsize \nu}}}}\!(x)\Delta t + o(\Delta t)
&
\!\!\!\ldots
&
\!\!
1+\partial_{{}_{\scriptstyle x_{{}_{\scriptsize
\xi_{{}_{\scriptsize \nu}}}}}}\!
f_{{}_{\scriptstyle \xi_{{}_{\scriptsize \nu}}}}\!\!(x)\Delta t + o(\Delta t)
\end{array}
\!\!\right|
=
\hfill
$
\\[3.55ex]
\mbox{}\hfill
$
= 1 +
{\rm div}{}^{\,\nu}_{\xi} f(x) \ \Delta t + o(\Delta t),
\hfill
$
\\[2ex]
where
$
{\rm div}{}^{\,\nu}_{\xi} f(x) =
\sum\limits_{j=1}^{\nu}
\partial_{{}_{\scriptstyle
x_{{}_{\scriptsize \xi_{{}_{\scriptsize j}}}}}}
f_{{}_{\scriptstyle \xi_{{}_{\scriptsize j}}}}(x)$
for all $x \in G.$

Therefore the first derivative of the scalar function
$
V^{\xi}\colon  {\rm U}(t_0) \to\R
$
is given by the formula
\\[2.25ex]
\mbox{}\hfill                                         
$
\displaystyle
{\sf D}V^{\xi}\colon t\to
\int\limits_{V^{\xi}_t}
{\rm div}{}^{\,\nu}_{\xi} f(x)
\,dx^{\nu}_{\xi}$
\ for all
$t \in {\rm U}(t_0^{}).
$
\hfill (1.1)
\\[-2.5ex]

\newpage

{\bf Lemma 1.1} [12, pp. 317 -- 319].
\vspace{0.25ex}
{\it
Suppose that $\lambda\!$-dimensional arcwise connected domain $G\subset \R^n,\,3\leq \nu \leq\lambda,$
has the homotopy group $\pi_{\nu-1}^{}(G)$ of rank $r,$
and for each sample $\xi=\xi_1^{}\ldots\,\xi_\nu^{}$ of dimension $\nu$ from the numbers $1,\ldots,n,$
\vspace{0.5ex}
there exists the single-valued scalar function
$\mu_{\xi}^{\nu}\in C^1(G)$ such that the
divergence
$
{\rm div}{}^{\,\nu}_{\xi}
M_{\xi}^{\nu}
$
of the single-valued vector function
}
\\[2ex]
\mbox{}\hfill              
$
M_{\xi}^{\nu}\colon x\to\ \mu_{\xi}^{\nu}(x)f(x)
$
\ for all \
$
x\in G
$
\hfill (1.2)
\\[2ex]
{\it
is $\!\nu\!$-constant sign on the domain $\!G.\!$
\vspace{0.5ex}
Then the following situation is impossible for $\!(\nu\! -\! 1)\!$-di\-men\-si\-o\-nal compact regular integral
manifolds of the ordinary  autonomous differential system {\rm (D)} with $f\in C^1(G)$
\vspace{0.35ex}
in any arcwise connected domain $G^{\;\!\prime}\subset G$ with homotopy group
$\pi_{\nu -1}^{}(G^{\;\!\prime})$
\vspace{0.5ex}
of rank $k\leq r\colon\!\!$
each of $(\nu - 1)\!$-di\-men\-si\-o\-nal compact regular integral manifolds $\Lambda^1,\ldots,\Lambda^k$
\vspace{0.5ex}
contains only one gap in its interior, and
$(\nu - 1)\!$-dimensional compact regular integral manifold
$\Lambda^{k + 1}$ contains all these $k$ gaps in its interior.
}
\vspace{0.5ex}

{\sl Proof}.
We use the following fact:
if the system {\rm (D)} with $f\in C^1(G)$ has a compact regular integral manifold of  dimension $\nu - 1$
in the arcwise connected domain $G,$ then
this manifold is $(\nu - 1)\!$-dimensional compact regular integral manifold in the domain $G$ of
the ordinary autonomous differential system with the vector field
$
M_{\xi}^{\nu} \colon x\to
\mu_{\xi}^{\nu}(x)f(x)
$
for all
$
x\in G,
$
where the scalar function $\mu_{\xi}^{\nu}\in C^1(G).$

Suppose the contrary, namely, let any of the
compact regular $(\nu - 1)\!$-dimensional integral manifold
$\Lambda^1,\ldots,\Lambda^{k}$ of system {\rm (D)}  with $f\in C^1(G)$
contain only one gap in its interior, and the
compact regular $(\nu - 1)\!$-dimensional integral manifold $\Lambda^{k+1}$ of system {\rm (D)}
contains all these $k$ gaps.
Then there exists a manifold $\Omega$ of dimension $\nu$
that is entirely inside the domain $G^{\;\!\prime}\subset G$
and is bounded by the manifolds $\Lambda^1,\ldots, \Lambda^k,\, \Lambda^{k+1}.$

Let $\R^{\nu}_{\xi},$  where the multiindex $\xi=\xi_1^{}\ldots\,\xi_n^{},$
be the subspace of the phase space $\R^n$ for which the natural projection
$\Omega^{\xi}$ of the manifold $\Omega$  has the dimension $\nu.$

Let us consider two logically possible cases:
\vspace{0.35ex}

1) the divergence ${\rm div}{}^{\,\nu}_{\xi}M^{\nu}_{\xi}$
on the domain $G^{\;\!\prime}$ is $\nu\!$-positive;
\vspace{0.75ex}

2) the divergence ${\rm div}{}^{\,\nu}_{\xi}M^{\nu}_{\xi}$
on the domain $G^{\;\!\prime}$ is $\nu\!$-negative.
\vspace{0.5ex}

In the first case, on the basis of the representation (1.1), we conclude that the
$\!\nu\!$-di\-men\-si\-o\-nal volume of the manifold $\Omega^{\;\!\xi}$ increases as $t \to {}+\infty.$
Since $\Lambda^1,\ldots,\Lambda^{k+1}$
are regular integral manifolds of dimensions ${\rm dim}\,\Lambda^j = \nu - 1, \, j=1,\ldots, k+1,\ \nu \geq 3,$
we arrive at a contradiction. Since the contradiction thus obtained takes place for any subspace $\R^{\nu}_{\xi}\subset\R^n$
for which the natural projection $\Omega^{\;\!\xi}$
of the manifold $\Omega$ has the dimension
${\rm dim}\,\Omega^{\;\!\xi} = \nu, \nu\geq 3,$
it follows that the assertion of the lemma holds for the case under consideration.

In the second case we arrive at a similar contradiction by substituting $t$ by ${}-t.$ \k
\vspace{0.5ex}

{\bf Theorem 1.1}
({\sl a test for the boundedness of the number of possible compact regular integral manifolds
of an ordinary autonomous differential system}) [10].
{\it
Suppose
that $\lambda\!$-di\-men\-si\-o\-nal arcwise connected domain $G^{\;\!\prime}\subset G$ with $3\leq \nu \leq\lambda$
has the homotopy group $\pi_{\nu-1}(G^{\;\!\prime})$ of rank $r$
\vspace{0.35ex}
and for each sample $\xi=\xi_1^{}\ldots\,\xi_\nu^{}$ of dimension $\nu$ from the numbers $1,\ldots,n$
\vspace{0.5ex}
there exists the single-valued scalar function
$\mu_{\xi}^{\nu}\in C^1(G^{\;\!\prime})$ such that the
divergence
$
{\rm div}{}^{\,\nu}_{\xi}
M_{\xi}^{\nu}
$
of the single-valued vector function {\rm (1.2)}
is $\!\nu\!$-constant sign on the domain $G^{\;\!\prime}.$
Then the ordinary autonomous differential system {\rm (D)} with $f\in C^1(G)$
has at most $r$ compact regular integral manifolds of dimension $\nu - 1$ in the domain  $G^{\;\!\prime}.$
}
\vspace{0.35ex}

{\sl Proof}.
For $r = 1,$ the desired assertion follows from Lemma 1.1.

Suppose that the assertion of  Theorem 1.1 holds for $r = k,$ i.e., for any
$\lambda\!$-di\-men\-si\-o\-nal arcwise connected domain in $\R^n$ with
homotopy group $\pi_{\nu - 1}^{}$ of rank $k.$

The following two logical cases are possible for a
$\lambda\!$-di\-men\-si\-o\-nal arcwise connected domain $G^{\;\!\prime}$ with
homotopy group $\pi_{\nu - 1}(G^{\;\!\prime})$ of rank $k+1\;\!\colon$

1) at least one gap is not inside a compact regular integral manifolds of dimension $\nu - 1;$

2) any gap is inside at least one compact regular integral manifolds of dimension $\nu - 1.$

In case 1), we split the domain $G^{\;\!\prime}$ into two $\lambda\!$-dimensional parts such that the above gap is
inside the boundary of one part, and the remaining gaps are inside the boundary of the other part.
Then, in the first part, there is no $(\nu - 1)\!$-dimensional compact regular integral manifolds of system (D)
with с $f\in C^1(G),$ and the other part contains at most $k - 1$
compact regular integral manifolds of this system (by assumption).

Therefore, in the first case, the domain $G^{\;\!\prime}$ contains at most $k - 1$
compact regular integral manifolds of dimension $\nu - 1.$

Let us consider case 2). Suppose that in the domain $G^{\;\!\prime}$ there are at least $k+1$
compact regular integral manifolds of dimension $\nu - 1.$
Then we have two logical possibilities:

a) no $(\nu - 1)\!$-dimensional compact regular integral manifold contains all gaps in its interior;

b) there exists a compact regular integral manifold of dimension $\nu - 1$
 that contains all gaps in its interior.

{\sl Case a}).
\vspace{0.2ex}
Suppose $\theta \geq 1$ gaps lie inside some $(\nu - 1)\!$-dimensional
compact regular integral manifold $\Lambda^i$ and
not exists $(\nu - 1)\!$-dimensional compact regular integral manifold $\Lambda^j$ such that
this manifold
also contain at least one of the  remaining  $k - \theta$ gaps.

We split the domain $G^{\;\!\prime}$ into two $\lambda\!$-dimensional parts such that the chosen
$(\nu - 1)\!$-di\-men\-si\-o\-nal compact regular integral manifold $\Lambda^i$
lies inside the boundary of one part, and no gap is outside this part.
Then, by the Lemma 1.1, a compact regular integral manifold of dimension $\nu - 1$
containing only these $\theta$ gaps can be neither inside this manifold $\Lambda^i$ nor outside it.
Taking account of the ranks of the homotopy groups $\pi_{\nu-1}^{}$ of the parts, we see
that the first part contains at most $\theta$ compact regular integral manifold of dimension $\nu - 1$
and the other at most $k-\theta$ ones.
Hence, in the domain $G^{\;\!\prime},$ there are at most $k$ compact regular integral manifold of dimension $\nu - 1.$

Thus case a) is impossible.
\vspace{0.25ex}

{\sl Case b}).
By the Lemma 1.1, all gaps cannot lie inside two compact regular integral manifold of dimension $\nu - 1.$
Then, by the same Lemma 1.1, inside the external compact regular integral manifold of dimension $\nu - 1$
there must be at least one gap not contained inside any
compact regular integral manifold of dimension $\nu - 1$ except for the external one.
Just as in case 1), we split the domain bounded by the external $(\nu - 1)\!$-dimensional compact regular integral manifold
into two $\lambda\!$-dimensional parts and see that the external
$(\nu - 1)\!$-dimensional compact regular integral manifold contains at most $k-1$ integral
$(\nu - 1)\!$-dimensional compact regular integral manifolds in its interior.
Hence, in the domain $G^{\;\!\prime},$ there are at most $k$ compact regular integral manifold of dimension $\nu - 1.$

Thus case b) is impossible.

The contradictions thus obtained mean that, in case 2), the domain $G^{\;\!\prime}$ contains at most $k$
compact regular integral manifold of dimension $\nu - 1$ of
\vspace{0.5ex}
system {\rm (D)} with $f\in C^1(G).$ \k

From Theorem 1.1, we get a
{\sl
test for the absence of nonisolated compact regular integral manifolds of
an ordinary autonomous differential system}.
\vspace{0.35ex}

{\bf Corollary 1.1.}\!\!
{\it
Under the assumptions of  Theorem {\rm 1.1}, in the domain $\!G^{\;\!\prime}\!\subset\! G\!$ system
{\rm (D)}  with $\!f\in C^1(G)\!$ has no
\vspace{0.75ex}
nonisolated $\!(\nu\! -\! 1)\!$-dimensional compact regular integral manifolds.
}

{\bf Example 1.1.}
The fifth-order  ordinary autonomous differential system
\\[2.25ex]
\mbox{}\hfill                                       
$
\begin{array}{c}
\dfrac{dx_1^{}}{dt} = {}-x_1^{} - x_2^{} + x_1^{}\;\!g(x)\equiv f_1^{}(x),
\qquad
\dfrac{dx_2^{}}{dt} = x_1^{} - x_2^{} + x_2^{}\;\!g(x)\equiv f_2^{}(x),
\\[3.5ex]
\dfrac{dx_3^{}}{dt} = {}-x_3^{} - x_4^{} + x_3^{}\;\!g(x)\equiv f_3^{}(x),
\qquad
\dfrac{dx_4^{}}{dt} = x_3^{} - x_4^{} + x_4^{}\;\!g(x)\equiv f_4^{}(x),
\\[3.5ex]
\dfrac{dx_5^{}}{dt} = {}-5x_5^{}\;\!g(x)\equiv f_5^{}(x),
\end{array}
$
\hfill (1.3)
\\[2.5ex]
where the scalar function
\vspace{0.25ex}
$
g\colon x\to x_1^2 + x_2^2 + x_3^2 + x_4^2$ for all $x\in \R^5,$
has three-dimensional compact regular integral manifold
\\[2ex]
\mbox{}\hfill               
$
\{x\colon g(x) = 1, \ x_5 = 0\}.
$
\hfill (1.4)
\\[2.15ex]
It can easily be checked that the derivative by virtue of system {\rm (1.3)}  is
\\[2.2ex]
\mbox{}\hfill
$
{\sf D}_t \bigl(g(x) - 1\bigr)_{\displaystyle |_{(1.3)}} =
2\bigl(g(x) - 1\bigr)g(x)
$
\ for all  $x\in \R^5.
\hfill
$
\\[2.3ex]
\indent
Let us prove its uniqueness (in the class of three-dimensional compact regular integral manifolds) with the help of Theorem 1.1.

The straight line
\vspace{0.5ex}
$
x_1^{} = x_2^{} = x_3^{} = x_4^{} = 0
$
is
a line of equilibria of system  (1.3).
The do\-ma\-in $\!G^{\;\!\prime}\!$ from $\R^5\backslash\{x\colon\! x_1^{} = x_2^{} = x_3^{} =x_4^{} = 0\}\!$
\vspace{0.5ex}
has the homotopy group $\!\pi_{3}^{}(G^{\;\!\prime})\!$ of rank 1.

Suppose
\\[1ex]
\mbox{}\hfill
$
\mu_{{}_{\scriptstyle 1234}}^{4}(x) = g^{{}-3}(x)
$
\
for all  $x\in G^{\;\!\prime},
\hfill
$
\\[2.5ex]
\mbox{}\hfill
$
\mu_{{}_{\scriptstyle 2345}}^{4}(x) =
\mu_{{}_{\scriptstyle 1345}}^{4}(x) =
\mu_{{}_{\scriptstyle 1245}}^{4}(x) =
\mu_{{}_{\scriptstyle 1235}}^{4}(x) = 1
$
\
for all  $x\in G^{\;\!\prime}.
\hfill
$
\\[2ex]
Then the divergences of the vector functions
\\[2.35ex]
\mbox{}\hfill
$
M^4_{{}_{\scriptstyle 1234}}(x)=
\mu^4_{{}_{\scriptstyle 1234}}(x)
\bigl(f_{1}(x),\ldots,f_{4}(x), 0\bigr)
$
\
for all $x\in G^{\;\!\prime},
\hfill
$
\\[2.5ex]
\mbox{}\hfill
$
M^4_{{}_{\scriptstyle 1235}}(x)=
\bigl(f_1(x),f_2(x),f_3(x), 0,f_5(x)\bigr)
$
\
for all $x\in G^{\;\!\prime},
\hfill
$
\\[2.5ex]
\mbox{}\hfill
$
M^4_{{}_{\scriptstyle 1245}}(x)=
\bigl(f_{1}(x),
f_{2}(x),0,f_{4}(x),f_{5}(x)\bigr)
$
\
for all
$x\in G^{\;\!\prime},
\hfill
$
\\[2.5ex]
\mbox{}\hfill
$
M^4_{{}_{\scriptstyle 1345}}(x)=
\bigl(f_{1}(x), 0,
f_{3}(x),f_{4}(x),f_{5}(x)\bigr)
$
\
for all $x\in G^{\;\!\prime},
\hfill
$
\\[2.5ex]
\mbox{}\hfill
$
M^4_{{}_{\scriptstyle 2345}}(x)=
\bigl(0,f_2(x),\ldots,f_5(x)\bigr)
$
\
for all $x\in G^{\;\!\prime}
\hfill
$
\\[2.75ex]
are $4\!$-definite functions on the domain $G^{\;\!\prime}\colon$
\\[2.75ex]
\mbox{}\hfill
$
{\rm div}^4_{{}_{\scriptstyle 1234}}\,
M^4_{{}_{\scriptstyle 1234}}(x) =
\dfrac{2}{(x_1^2 + x_2^2 + x_3^2 + x_4^2)^{3}}  >  0\;\!;
\hfill
$
\\[2.75ex]
\mbox{}\hfill
$
{\rm div}^4_{{}_{\scriptstyle 1235}}\,
M^4_{{}_{\scriptstyle 1235}}(x) = {} -(3 + 2x_4^2) < 0\;\!;
\qquad
{\rm div}^4_{{}_{\scriptstyle 1245}}\,
M^4_{{}_{\scriptstyle 1245}}(x) = {}-(3 + 2x_3^2) < 0\;\!;
\hfill
$
\\[2.75ex]
\mbox{}\hfill
$
{\rm div}^4_{{}_{\scriptstyle 1345}}\,
M^4_{{}_{\scriptstyle 1345}}(x) ={} -(3 + 2x_2^2) < 0\;\!;
\qquad
{\rm div}^4_{{}_{\scriptstyle 2345}}\,
M^4_{{}_{\scriptstyle 2345}}(x) = {}-(3 + 2x_1^2) < 0.
\hfill
$
\\[3ex]
\indent
By Theorem 1.1, the manifold (1.4) is the unique three-dimensional compact regular integral manifold
of the ordinary autonomous differential system (1.3).
\\[4.75ex]
\centerline{
{\bf 1.2. Autonomous total differential system}
}
\\[2.25ex]
\indent
{\bf 1.2.1.\!\! Test of the boundedness of the number of compact regular integral ma\-ni\-folds.}\!
Autonomous system of total differential equations {\rm (TD)} induces $m$ ordinary autonomous
differential systems of $n\!$-th order
\\[2ex]
\mbox{}\hfill                                        
$
dx = X^j(x)\,dt_{j}^{}\;\!,
\quad
j = 1,\ldots, m,
$
\hfill {\rm (Dj)}
\\[2.25ex]
where $X^j(x) ={\rm colon}\bigl(X_{1j}^{}(x), \ldots, X_{nj}^{}(x)\bigr)$ for all $x\in G.$
\vspace{1ex}

{\bf Lemma 1.2.}
{\it
If the system {\rm (TD)} has a regular integral manifold, then each of the systems {\rm (Dj)}, $j=1,\ldots,  m,$
has the same regular integral manifold. In addition, the compactness of regular integral manifolds is preserved.
}
\vspace{0.35ex}

We can prove Lemma 1.2
\vspace{0.75ex}
by the method of fixing $m-1$ independent variables $t_1^{},\ldots,t_m^{}.$

Using  Lemma 1.2, from Theorem 1.1, we get a
{\sl test of the boundedness of the number of compact regular integral manifolds for
\vspace{0.5ex}
an autonomous total differential system} [10].

{\bf Theorem 1.2.}
{\it
Suppose
that $\lambda\!$-di\-men\-si\-o\-nal arcwise connected domain $G^{\;\!\prime}\subset G$ with $3\leq \nu \leq\lambda$
has the homotopy group $\pi_{\nu-1}(G^{\;\!\prime})$ of rank $r$
\vspace{0.35ex}
and for each sample $\xi=\xi_1^{}\ldots\,\xi_\nu^{}$ of dimension $\nu$ from the numbers $1,\ldots,n$
\vspace{0.5ex}
there exists the single-valued scalar function
$\mu_{\xi j}^\nu \in C^1(G^{\;\!\prime}),\ j\in\{1,\ldots,m\},$
\vspace{0.35ex}
such that the
divergence
$
{\rm div}{}^{\,\nu}_{\xi}\;
M_{\xi j}^{\nu}
$
of the single-valued vector function
$
M_{\xi  j}^{\nu}\colon x \to
\mu_{\xi  j}^{\nu}(x)X^j(x)$ for all $x \in G^{\;\!\prime}$
\vspace{0.35ex}
is $\!\nu\!$-constant sign on the domain $G^{\;\!\prime}.$
Then the autonomous total differential  system {\rm (TD)} with $X\in C^1(G)$
\vspace{0.25ex}
has at most $r$ compact regular integral manifolds of dimension $\nu - 1$ in the domain  $G^{\;\!\prime}.$
}
\vspace{0.35ex}

Using the estimates of Theorem 1.2 and successively considering the ordinary autonomous differential systems
{\rm (D1)}, \ldots, {\rm (D}m), we can estimate the total possible number of  compact regular integral manifolds of system (TD).
\vspace{0.35ex}

From Theorem 1.2, we obtain
{\sl a test of the absence of nonisolated compact regular integral manifolds for
an autonomous total differential system}.
\vspace{0.35ex}

{\bf Corollary 1.2.}
{\it
If there exists $j \in \{1,\ldots, m\}$
such that the system {\rm (Dj)}
satisfies the as\-sum\-p\-ti\-ons of Theorem {\rm 1.2}, then the system {\rm (TD)} with $X\in C^1(G)$
has no nonisolated compact regular integral manifolds of dimension $\nu - 1$ in the domain  $G^{\;\!\prime}\subset G.$
}
\vspace{0.5ex}

{\bf Example 1.2.}
For the autonomous system of total differential equations
\\[2.25ex]
\mbox{}\hfill                                         
$
\begin{array}{l}
dx_1^{} = \bigl({}-x_1^{} - x_2^{} + x_1^{}g(x)\bigr)\;\!dt_1^{}
+ \bigl({}-x_1^{} + x_4^{} + x_1^{}g(x)\bigr)\;\! dt_2^{},
\\[2.1ex]
dx_2^{} = \bigl(x_1^{} - x_2^{} + x_2^{}g(x)\bigr)\;\!dt_1^{} +
\bigl({}-x_2^{} + x_3^{} + x_2^{}g(x)\bigr)\;\! dt_2^{},
\\[2.1ex]
dx_3^{} =  \bigl({}-x_3^{} - x_4^{} + x_3^{}g(x)\bigr)\;\! dt_1^{}  +
\bigl({}-x_2^{} - x_3^{} + x_3^{}g(x)\bigr)\;\!dt_2^{},
\\[2.1ex]
dx_4^{} = \bigl(x_3^{} - x_4^{} + x_4^{}g(x)\bigr)\;\! dt_1^{} +
\bigl({}-x_1^{} - x_4^{} + x_4^{}g(x)\bigr)\;\! dt_2^{},
\end{array}
$
\hfill (1.5)
\\[2.5ex]
where
$
g\colon x\to x_1^2 + x_2^2 + x_3^2 + x_4^2$ for all $x\in\R^4,$
the sphere
$
S^3 = \{x\colon g(x)= 1\}
$
is a three-dimensional compact regular integral manifold.
We see that the differential by virtue of the total differential system {\rm (1.5)}  is
\\[2ex]
\mbox{}\hfill
$
d\bigl(g(x) - 1\bigr)_{\displaystyle |_{(1.5)} }  =
2\bigl(g(x) - 1\bigr) g(x)(dt_1 + dt_2)
$
\ for all
$(t,x)\in \R^6.
\hfill
$
\\[2.5ex]
\indent
Let us consider the ordinary differential system
\\[2.5ex]
\mbox{}\hfill                                         
$
\begin{array}{l}
dx_1^{} = \bigl({}-x_1^{} - x_2^{} + x_1^{}g(x)\bigr)\;\!dt_1^{},
\qquad
dx_2^{} = \bigl(x_1^{} - x_2^{} + x_2^{}g(x)\bigr)\;\!dt_1^{},
\\[2.75ex]
dx_3^{} =  \bigl({}-x_3^{} - x_4^{} + x_3^{}g(x)\bigr)\;\! dt_1^{},
\qquad
dx_4^{} = \bigl(x_3^{} - x_4^{} + x_4^{}g(x)\bigr)\;\! dt_1^{}
\end{array}
$
\hfill (1.6)
\\[2.5ex]
induced by the total differential system (1.5).
\vspace{0.35ex}

Let
$
\mu_{{}_{\scriptstyle 1234}}^{4}\colon x \to g^{{}-3}(x)$ for all $x \in G^{\;\!\prime},$
where a domain $G^{\;\!\prime}\subset \R^4\backslash\{(0,0,0,0)\}.
$
Then, for the ordinary autonomous differential system {\rm (1.6)},
the divergence of the vector function $M_{{}_{\scriptstyle 1234}}^{4}$
is 4-positive sign on the domain $G^{\;\!\prime}\colon$
\\[2.2ex]
\mbox{}\hfill
$
{\rm div}\,M^{4}_{{}_{\scriptstyle 1234}}(x) =
2  g^{{}-3}(x)  >  0
$
\ for all
$
x\in G^{\;\!\prime}.
\hfill
$
\\[2.65ex]
\indent
Since the rank of the homotopy group $\!\pi_{3}^{}(G^{\;\!\prime})\!$ is $1$ we see that
the sphere $S^3\!=\!\{x\colon\! g(x)\!=\! 1\}$ is the unique three-dimensional compact regular integral manifold of
the ordinary auto\-no\-mo\-us differential system {\rm (1.6)} (by Theorem 1.1).

By Theorem 1.2,
the sphere $S^3=\{x\colon g(x)= 1\}$  is
the unique three-dimensional compact regular integral manifold of the total differential system (1.5).
\\[3.25ex]
\indent
{\bf 1.2.2.\! Tests of the absence of compact regular orbits}.
Consider the completely solvable autonomous system of total differential equations
{\rm (TD)} with $X\in C^k(G), \ 1\leq k\leq \infty.$
Let us remember that
the Frobenius conditions [13] for completely solvability
of sys\-tem (TD) we may represent via Poisson brackets as the identities
\\[2ex]
\mbox{}\hfill       
$
[ {\mathfrak X}_j^{}(x), {\mathfrak X}_l^{}(x)] = {\mathfrak O}
$
\ for all $x\in G,
\ \ j=1,\ldots,m, \  \ l=1,\ldots,m,
$
\hfill (1.7)
\\[2ex]
where ${\frak O}$ is the null operator,
the linear differential ope\-ra\-tors
\\[1.5ex]
\mbox{}\hfill
$
\displaystyle
{\mathfrak X}_j^{}(x)=\sum\limits_{i=1}^n X_{ij}^{}(x)\,\partial_{x_i^{}}^{}
$
\ for all $x\in G,
\quad
j=1,\ldots,m,
\hfill
$
\\[1.5ex]
induced by system (TD) are operators of differentiation by virtue of sys\-tem (TD).
\vspace{0.75ex}

{\bf Theorem 1.3} ({\sl the first test of the absence of compact regular orbits for
an autonomous total differential system}) [11].
\vspace{0.25ex}
{\it
Suppose there exists a single-valued scalar function
$N\in C^1(G^{\;\!\prime}),\, G^{\;\!\prime}\subset G,$ such that
\\[1.75ex]
\mbox{}\hfill                                         
$
{\frak X}_{{}_{\scriptstyle j}}N(x) =H_{j}^{}(x)
$
\ for all $
x\in G^{\;\!\prime},
\quad
j =1,\ldots, m,
$
\hfill {\rm(1.8)}
\\[1.75ex]
and the function
\vspace{0.35ex}
$
H_{k}^{}\colon  G^{\;\!\prime}\to \R
$
is definite for at last one $k\in\{1,\ldots, m\}.$ Then the
completely solvable system {\rm (TD)} with $X\in C^{\infty}(G)$ has no
compact regular orbits in the domain $G^{\;\!\prime}.$
}
\vspace{0.5ex}

{\sl Proof}.
Suppose that in the domain $G^{\;\!\prime},$ there is a compact regular orbit of the completely solvable
total differential system {\rm(TD)} with $X\in C^{\infty}(G)$
\vspace{0.25ex}
corresponding to a periodic solution
$
x\colon t \to x(t)$  for all
$
t \in  D
$
with period
\vspace{0.35ex}
$T = (T_1,\ldots, T_m)$ and
initial condition
$
x_{0}^{} =x(t_0^{}),\  t_0^{} =
 (t_{01}^{},\ldots, t_{0m}^{}), \  t_0^{}\in D, \
x_0^{} \in G^{\;\!\prime}.
$
Then
\vspace{0.5ex}
$
T_{j}^{} \ne 0,\  j =1,\ldots,m,
$
and, by Theorem 5.4 in [13, c. 31],
this solution can be continued to the entire space $\R^m.$
\vspace{0.5ex}

Let the function
\vspace{0.35ex}
$H_{k}^{}\colon G^{\;\!\prime}\to {\R}$ be positive definite (negative definite).
Then from the system of identities (1.8), we find that
if
\vspace{0.75ex}
$
t_{j}^{} = t_{0j}^{},\
j=1,\ldots,m,\ j \ne k,\,\
t_{0k}^{} \leq t_{k}^{} \leq
t_{0k}^{} + T_{k}^{},
$
then the function
$
v\colon t_{k}^{} \to N(x(t))
$
\vspace{0.5ex}
is strictly increasing (strictly decreasing).
Therefore,
$
v(t_{ 0k}^{}) < v(t_{ 0k}^{} + T_{k}^{})\
 \
\bigl(v(t_{0k}^{}) >
v(t_{0k}^{} + T_{k}^{})\bigr),
$
\vspace{0.5ex}
which, together with the single-valuedness of $N,$ contradicts the periodicity of solution
$x\colon t\to x(t)$ for all $t \in D.$
\vspace{0.25ex}

This contradiction completes the proof of Theorem 1.3. \k
\vspace{0.35ex}

The first test (Theorem 1.3) is coordinated with the Poincar\'{e} test
and his generalized version (Theorem 0.1)
when the generalized Lyapunov function used to determine
the absence of closed trajectories for the ordinary autonomous
\vspace{0.35ex}
differential system of the second order (0.1).

We introduce the linear differential form
\\[2ex]
\mbox{}\hfill                                          
$
\displaystyle
\omega(x) = \sum\limits_{i = 1}^{n}\,
w_i^{}(x)\,dx_i^{}$
\ for all $x\in G^{\;\!\prime}
$
\hfill  (1.9)
\\[2ex]
that is exact on the domain $G^{\;\!\prime}\subset G$ and has
continuously differentiable on the domain $G^{\;\!\prime}$
coefficients $w_{i}^{}\colon G^{\;\!\prime}\to \R,\, i=1,\ldots, n,$ i.e.,
$\omega\in C^1(G^{\;\!\prime}).$

Sinse the 1-form (1.9) is exact on the domain $G^{\;\!\prime},$ we see that
there exists a single-valued function $N\colon G^{\;\!\prime}\to {\R}$
shuch that the total differential
\\[1.5ex]
\mbox{}\hfill
$
dN(x) = \omega(x)$
\ for all $x \in G^{\;\!\prime}.
\hfill
$
\\[0.75ex]
Therefore,
\\[1ex]
\mbox{}\hfill                                          
$
\displaystyle
{\frak X}_{j}^{}N(x) =
\sum\limits_{i=1}^{n}\,
w_i^{}(x) X_{ij}^{}(x)
$
\ for all $x \in G^{\;\!\prime},
\quad
j = 1,\ldots, m.
$
\hfill  (1.10)
\\[1.75ex]
\indent
Using Theorem 1.3, we obtain
{\sl the second test of the absence of compact regular orbits for an autonomous total differential system} [11].
\vspace{0.5ex}

{\bf Theorem 1.4.}
{\it
\vspace{0.35ex}
If there exists on the domain $G^{\;\!\prime}\subset G$
a linear differential form {\rm(1.8)} such that the sum
$
\sum\limits_{i=1}^{n}
w_i^{}(x)
X_{ij}^{}(x)
$
\vspace{0.35ex}
is definite on the domain $G^{\;\!\prime}$ for at least one  $j\in\{1,\ldots, m\},$
then the system {\rm (TD)} with $X\in C^{\infty}(G)$ has no
compact regular orbits in the domain $G^{\;\!\prime}.$
}
\vspace{0.5ex}

Note that Theorem 1.4 (unlike Theorem 1.3) does not require knowledge of the function $N$ and
only assumes that the derivative (1.10) of $N$ by virtue of one of the ordinary autonomous differential systems
{\rm (Dj)} induced by the completely solvable system {\rm (TD)} with $X\in C^{\infty}(G)$
is definite on the domain $G^{\;\!\prime}.$

Let us consider a linear homogeneous autonomous system of total differential equations
\\[2ex]
\mbox{}\hfill       
$
dx=A(x)\, dt,
$
\hfill (1.11)
\\[2.15ex]
where
\vspace{0.75ex}
$t\in\R^m,\ x\in\R^n,\ dx={\rm colon}\, (dx_1,\ldots,dx_n),\
dt={\rm colon}\, (dt_1,\ldots,dt_m),\ m<n,$
the entries of
$n\times m$
matrix
$A(x)=\bigl|\!\bigl| a_{ij}^{}(x) \bigl|\!\bigl|$ for all $x\in\R^n$
are linear homogeneous functions
$
a_{ij}^{} \colon x\to
\sum\limits_{\tau=1}^{n} a_{ij\tau}^{}x_\tau^{}$
for all
$x\in\R^n,
\ a_{ij\tau}^{}\in\R,\
i=1,\ldots,n,\  j=1,\ldots,m,
\ \tau =1,\ldots,n.
$
\vspace{0.5ex}

The Frobenius conditions for completely solvability of system (1.11)  are equivalent
\\[1.5ex]
\mbox{}\hfill
$
A_j^{}A_l^{}=A_l^{}A_j^{},
\quad
 j=1,\ldots,m,\  \ l=1,\ldots,m,
\hfill
$
\\[1.5ex]
where the $\!n\times n\!$-matrices
\vspace{0.75ex}
$\!A_j\!=\!\bigl|\!\bigl| a_{\tau  j i }^{} \bigl|\!\bigl|,\, j\!=\!1,\ldots,m\!$
(\!$\tau\!$ is a row number, $\!i$ is a column number).

By Theorem 1.1 in [20, p. 30] (or Lemma 7.5.1 in [21, p. 249])
and  Theorem 1.3, we have
{\sl  a test of the absence of compact regular orbits for a linear homogeneous autonomous total differential system} [11].
\vspace{0.5ex}

{\bf Theorem 1.5.}
{\it
\vspace{0.15ex}
If at least one of the matrices $A_{j}^{}$
has eigenvalues $\lambda^{j}_{1},\ldots, \lambda^{j}_{n}$
such that
$
\lambda^{j}_{i} + \lambda^{j}_{\tau}\ne 0, \
i=1,\ldots,n, \, \tau=1,\ldots,n,\,
j\in\{1,\ldots, m\},
$
\vspace{0.35ex}
then the completely solvable linear homogeneous autonomous total differential system {\rm (1.11)}
has no compact regular orbits.
}
\\[5.75ex]
\centerline{
\large\bf
2.\! The boundedness of the number of compact integral hypersurfaces
}
\\[2.25ex]
\centerline{
{\bf  2.1. System of exterior differential equations}
}
\\[1.75ex]
\indent
{\bf Theorem 2.1}
({\sl a test for the boundedness of the number of possible
compact integral hypersurfaces of a system of exterior differential equations}).
\vspace{0.35ex}
{\it
Let a domain $G^{\;\!\prime}\subset G$  have the homotopy group $\pi_{n-1}(G^{\;\!\prime})$ of rank $r.$
Suppose that there exist an
\vspace{0.5ex}
$(n-2)\!$-form $\alpha\in C^2(G^{\;\!\prime})$ and
$(n-p_{{}_{\scriptstyle j}}-1)\!$-forms $\gamma_{j}\in C^1(G^{\;\!\prime}),\ j=1,\ldots, m,$
such that the exterior differentials of sum
\\[2ex]
\mbox{}\hfill                                      
$
\displaystyle
d\Bigl(d\alpha(x)_{\displaystyle |_{\scriptstyle {\rm (ED)} }} +
\sum\limits_{j=1}^{m}\,
\zeta_{_{\scriptstyle j}}(x) \wedge
\gamma_{_{\scriptstyle j}}(x) \Bigr) =
B(x)\,dx_1 \wedge \ldots \wedge
dx_n,
$
\hfill {\rm(2.1)}
\\[2.25ex]
where the function $B\colon G^{\;\!\prime}\to\R$ is constant sign.
\vspace{0.35ex}
Then in the domain $G^{\;\!\prime}$ the system
{\rm (ED)} with $\zeta_{j}\in C^{\infty}(G),\ j=1,\ldots, m,$
has at most r compact integral hypersurfaces.
}
\vspace{0.75ex}

{\sl Proof} of Theorem 2.1 [9]
is similar to that one in Theorem 1.1 when $\nu=n.$
Also, we use the following
\vspace{0.35ex}

{\bf Lemma 2.1.}
{\it
Let the assumptions of Theorem {\rm 2.1} be valid.
Then the following arrangement is impossible for compact integral hypersurfaces of the system of external differential equations
{\rm (ED)} with $\zeta_{j}^{}\in C^{\infty}(G),\ j=1,\ldots, m,$ in
any subdomain $\Omega$ of the domain $G^{\;\!\prime}\subset G$
with homotopy group $\pi_{n-1}^{} (\Omega)$ of rank $s,\ s\leq r\colon$
any of $s$ gaps is surrounded by its own compact integral hypersurface
$\partial\Sigma_{1}^{},\ldots, \partial\Sigma_{s}^{},$ and a compact integral hypersurface
$\partial\Sigma_{s+1}^{}$ sur\-ro\-unds these $s$ gaps{\rm;}
moreover, the hypersurfaces $\partial\Sigma_{1}^{},\ldots, \partial\Sigma_{s}^{}$ are disjoint, do not sur\-ro\-und each other,
and are entirely surrounded by the hypersurface $\partial\Sigma_{s+1}^{}.$
}
\vspace{0.5ex}

{\sl Proof}.
We prove Lemma 2.1 by contradiction. Assume that the situation described in Lemma 2.1 takes place.
\vspace{0.25ex}
Then by $\Sigma$  we denote the domain bounded by the hypersurface
\linebreak
$\partial\Sigma = \bigcup\limits^{s+1}_{\tau=1}\partial\Sigma_{{}_{\scriptstyle \tau}}.$
\vspace{0.35ex}
Taking into account the Stokes formula for an oriented manifold with boundary and the identity (2.1), we obtain
\\[2ex]
\mbox{}\hfill                                              
$
\displaystyle
\int\limits_{\partial\Sigma}
\Bigl( d\alpha(x)_{\displaystyle |_{\scriptstyle {\rm (ED)}}}
+ \sum\limits_{j=1}^{m}\,
\zeta_{_{\scriptstyle j}}(x) \wedge
\gamma_{_{\scriptstyle j}}(x) \Bigr) =
\hfill
$
\\[-2ex]
\mbox{}\hfill (2.2)
\\[2.75ex]
\mbox{}\hfill
$
\displaystyle
= ({}-1)^n \int\limits_{\Sigma}
d\Bigl( d\alpha(x)
_{\displaystyle |_{\scriptstyle {\rm (ED)}}}
\,  +\,
\sum\limits_{j=1}^{m}\,
\zeta_{_{\scriptstyle j}}(x) \wedge
\gamma_{_{\scriptstyle j}}(x) \Bigr)
\, =\,   ({}-1)^n \int\limits_{\Sigma} B(x)\, d\Sigma.
\hfill
$
\\[1.75ex]
\indent
Using the Poincare theorem [15, p. 111] to the effect that
$d(d\omega(x))=0$ for all $x\in G,$ where the differential form $\omega\in C^2(G),$
in view of the system of exterior differential equations (ED) we have
\\[1.5ex]
\mbox{}\hfill                                            
$
\displaystyle
\int\limits_{\partial\Sigma}
\Bigl( d\alpha(x)_{\displaystyle |_{\scriptstyle {\rm (ED)}}}
\, +\,
\sum\limits_{j=1}^{m}\,
\zeta_{_{\scriptstyle j}}(x) \wedge \gamma_{_{\scriptstyle j}}(x)
\Bigr)_
{\displaystyle |_{\scriptstyle {\rm (ED)}}}
 =
\hfill
$
\\[-2ex]
\mbox{}\hfill (2.3)
\\[2.75ex]
\mbox{}\hfill
$
\displaystyle
=({}-1)^n \int\limits_\Sigma
d\Bigl( d\alpha(x)_
{\displaystyle |_{\scriptstyle {\rm (ED)}}}\Bigr)
\ + \
\int\limits_{\partial\Sigma}\
\sum\limits_{j=1}^{m}\;\!\bigl(
\zeta_{_{\scriptstyle j}}(x) \wedge
\gamma_{_{\scriptstyle j}}(x)\bigr)
_{\displaystyle |_{\scriptstyle {\rm (ED)}}}
 = 0.
\hfill
$
\\[2ex]
\indent
By virtue of (2.2), the equality (2.3) is impossible, since the integrand of the multiple integral occurring on the right-hand side in the chain of relations (2.2) is a constant sign function on the domain $G^{\;\!\prime}$ and $\Sigma\subset\Omega\subset G^{\;\!\prime}.$
The obtained contradiction completes the proof of Lemma 2.1. \k

{\bf Example 2.1.}
Consider the system of exterior differential equations
\\[1.75ex]
\mbox{}\hfill                                       
$
\zeta_j^{}(x)=0,
\quad
j=1,\ldots, 4,
$
\hfill (2.4)
\\[1.75ex]
where the differential forms
\\[2ex]
\mbox{}\hfill
$
\zeta_1^{}(x)=
(x_1^2 + x_2^2 + x_3^2 + x_4^2)
\bigl(x_2^2\,dx_1 + (x_3^{} +x_4^2)\,dx_2\bigr) +
(x_1^{} + x_2^{} - x_3^{} + x_4^{} + x_1^2 - 2x_1^{}x_2^{} +3x_4^2)\,dx_3^{} \ +
\hfill
$
\\[2ex]
\mbox{}\hfill
$
+\ (2x_2^{} - x_3^{} + 5x_4^{} + 3x_1^2 - x_2^2 + 2x_3^2 - 5x_4^2)\,dx_4^{}
$
\ for all
$
x\in\R^4,
\hfill
$
\\[2.25ex]
\mbox{}\hfill
$
\zeta_2^{}(x)=
dx_1^{} + x_1^{}({}-2x_2^{} + x_1^2)\,dx_2^{} + (x_1^2 + x_3^2)\,dx_3^{} +(x_2^2 + x_4^2)\,dx_4^{}
$
\ for all
$x\in\R^4,
\hfill
$
\\[2.25ex]
\mbox{}\hfill
$
\zeta_3^{}(x) =
\bigl({}-1 + 2x_1^{} + (x_1^2 + x_2^2 + x_3^2 + x_4^2)(1 - x_2^2)\bigr)\;\!dx_1^{}\ +
\hfill
$
\\[2ex]
\mbox{}\hfill
$
+\
\bigl(5 + 2x_2^{} - (x_1^2 + x_2^2 + x_3^2 + x_4^2)(5 + x_3^{} + x_4^2)\bigr)\!\;dx_2^{} \,+\,
({}-x_1^{} - x_2^{} + 3x_3^{} - x_4^{} - x_1^2 + 2x_1^{}x_2^{} - 3x_4^2)\, dx_3^{} \ +
\hfill
$
\\[2ex]
\mbox{}\hfill
$
+\
({}-2x_2^{} + x_3^{} - 3x_4^{} - 3x_1^2 +x_2^2 -2x_3^2  + 5x_4^2)\,dx_4^{}
$
\ for all
$x\in\R^4,
\hfill
$
\\[2.25ex]
\mbox{}\hfill
$
\zeta_4^{}(x)=
x_1^2\;dx_1^{} \wedge dx_4^{} + x_2^2\; dx_2 \wedge dx_3^{}
$
\ for all $x\in\R^4.
\hfill
$
\\[2ex]
\indent
The system (2.4) such that the differential 2-forms
\\[2ex]
\mbox{}\hfill
$
\alpha(x)  =  x_1^{}\, dx_3^{}\wedge dx_4^{}$ for all $x\in\R^4,
\hfill
$
\\[2.5ex]
\mbox{}\hfill
$
\gamma_{1}^{}(x) =
\dfrac{1}{x_1^2+x_2^2 + x_3^2 + x_4^2}\ dx_3^{}\wedge dx_4^{}
$
\ for all
$x\in G^{\;\!\prime},
\quad
\gamma_2^{}(x) =
\gamma_3^{}(x) = 0
$
\ for all $x\in \R^4,
\hfill
$
\\[2.5ex]
and the differential 1-form
\\[0.5ex]
\mbox{}\hfill
$
\gamma_4^{}(x) = 0$ for all $x\in \R^4
\hfill
$
\\[2ex]
on the domain $G^{\;\!\prime}=\R^4\backslash\{(0,0,0,0)\}$
satisfy the relations
\\[2ex]
\mbox{}\hfill
$
d\alpha(x)_{\displaystyle |_{\scriptstyle(2.4)}} =
x_1^{}(2x_2^{} -  x_1^2)\,dx_2^{} \wedge dx_3^{} \wedge dx_4^{},
\hfill
$
\\[2.75ex]
\mbox{}\hfill
$
\zeta_1^{}(x) \wedge
\gamma_1^{}(x) =
x_2^2\;dx_1^{} \wedge dx_3^{} \wedge dx_4^{} +
(x_3^{} + x_4^2)\,dx_2^{}\wedge dx_3^{} \wedge dx_4^{},
\hfill
$
\\[3ex]
\mbox{}\hfill
$
\zeta_2^{}(x)\wedge
\gamma_2^{}(x) =
\zeta_3^{}(x) \wedge
\gamma_3^{}(x) =
\zeta_4^{}(x) \wedge
\gamma_4^{}(x) = 0,
\hfill
$
\\[2ex]
\mbox{}\hfill
$
\displaystyle
d\Bigl(d\alpha(x)_{\displaystyle |_{\scriptstyle(2.4)}} +
\sum\limits_{j=1}^{4}\,
\zeta_j^{}(x) \wedge
\gamma_j^{}(x) \Bigr) =
{}-3x_1^2\;dx_1^{} \wedge dx_2^{} \wedge dx_3^{} \wedge dx_4^{}.
\hfill
$
\\[1.75ex]
\indent
By Theorem 2.1, in the domain $G^{\;\!\prime}=\R^4\backslash\{(0,0,0,0)\}$
\vspace{0.35ex}
with homotopy group $\pi_{{}_{\scriptstyle 3}}(G^{\;\!\prime})$ of rank 1 the system (2.4)
can have at most one compact integral hypersurface.

Since
\\[0.75ex]
\mbox{}\hfill
$
dw(x) = \bigl(\zeta_1(x) +
\zeta_3(x)\bigr)_{\displaystyle |_{\scriptstyle w(x)=0}},
\hfill
$
\\[2ex]
where
\vspace{0.35ex}
$
w\colon x\to x_1^2 + x_2^2 + x_3^2 +
x_4^2 - 1
$
for all
$
x\in\R^4,
$
we see that the sphere
$
S^3= \{x\colon w(x)=0\}
$
is this single compact integral hypersurface of system (2.4).
\vspace{0.5ex}

The possibility of the fact that some gap or a set of gaps of the domain $G^{\;\!\prime}$
is not surrounded by a compact integral hypersurface of the system of exterior differential equations {\rm (ED)}
with $\zeta_{j}\in C^{\infty}(G),\ j=1,\ldots, m,$
can be considered on the basis of the notion of invariance of a differential form on the domain $G^{\;\!\prime}$
with respect to system (ED).

A differential $(n-2)\!$-form $\theta$ is said to be {\it invariant} on the domain $G^{\;\!\prime}\subset G$
under the system of exterior differential equations (ED) if it is closed on any $(n-2)\!$-dimensional manifold of this system, i.e.,
\\[1.5ex]
\mbox{}\hfill
$
d\;\!\theta(x)_{\displaystyle |_{\scriptstyle {\rm (ED)}}} =  0
$
\ for all
$
x \in G^{\;\!\prime},
\hfill
$
\\[1.75ex]
or, in other words, there exist $(n-p_{{}_{\scriptstyle j}}-1)\!$-forms
\vspace{0.5ex}
$\eta_{{}_{\scriptstyle j}},\, j = 1,\ldots, m,$
with continuously differentiable coefficients on the domain $G^{\;\!\prime}$
and such that the exterior differential
\\[1.75ex]
\mbox{}\hfill
$
\displaystyle
d\;\!\theta(x) = \sum\limits_{j=1}^{m}\,
\zeta_{{}_{\scriptstyle j}}(x) \wedge
\eta_{{}_{\scriptstyle j}}(x)
$
\ for all $x\in G^{\;\!\prime}.
\hfill
$
\\[1.75ex]
\indent
The {\it index of the gap} $\Theta$ of the domain $G^{\;\!\prime}$
\vspace{0.55ex}
with respect to the closed differential $(n-1)\!$-form $\delta$ on the domain $G^{\;\!\prime}$
is called the number
$
{\rm ind}_{{}_{\scriptstyle \delta}}\Theta =
\int\limits_S\delta,
$
where $S$ is a manifold homeomorphic to a hypersphere lying in the domain $G^{\;\!\prime}$
\vspace{0.25ex}
such that the gap $\Theta$ is the only obstruction to the continuous contraction of this manifold into a point.
\vspace{0.75ex}

{\bf Theorem 2.2} [9].
\vspace{0.25ex}
{\it
Let a domain $G^{\;\!\prime}\subset G$ have the homotopy group $\pi_{n-1}(G^{\;\!\prime})$ of rank $r,$
and let there exist differential $(n-2)\!$-forms $\alpha\in C^2(G^{\;\!\prime}),\,\theta\in C^2(G^{\;\!\prime}),$ and
$(n-p_{{}_{\scriptstyle j}}-1)\!$-forms $\gamma_j\in C^1(G^{\;\!\prime}),\, j= 1,\ldots,m,$
such that the $(n-2)\!$-form $\theta$ is invariant on the domain $G^{\;\!\prime}$
\vspace{0.25ex}
with respect to the system of exterior differential equations {\rm (ED)} with $\zeta_{j}\in C^{\infty}(G),\ j=1,\ldots, m,$
and the exterior differential of sum}
\\[2ex]
\mbox{}\hfill                                            
$
\displaystyle
d\Bigl(d\alpha(x)_
{\displaystyle |_{\scriptstyle {\rm (ED)}}}\! + d\theta(x) +
\sum\limits_{j=1}^{m}\,
\zeta_{_{\scriptstyle j}}(x) \wedge
\gamma_{_{\scriptstyle j}}(x) \Bigr) =
B(x)\,dx_1^{} \wedge \ldots \wedge dx_n^{}
$
\ for all
$
x\in G^{\;\!\prime},
$
\hfill (2.5)
\\[2ex]
{\it
where the function $B\colon G^{\;\!\prime}\to\R$ is constant sign. Then  the following assertions are valid{\rm:}
\vspace{0.35ex}

{\rm 1)}
\vspace{0.25ex}
the system {\rm (ED)} can have at most $r$  compact integral hypersurfaces in the domain $G^{\;\!\prime};$

{\rm 2)}
any set of gaps of the domain $G^{\;\!\prime}$ surrounded by a compact integral hypersurface of system {\rm (ED)}
has the zero total index with respect to the $(n-1)\!$-form $d\theta.$
}
\vspace{0.35ex}

{\sl Proof.}
The first assertion readily follows from Theorem 2.1. Since
on the domain $G^{\;\!\prime}$ the exterior differential
\\[2ex]
\mbox{}\hfill
$
\displaystyle
d\Bigl(
d\alpha(x)_{\displaystyle |_{\scriptstyle {\rm (ED)}}}
+ d\theta(x) +
\sum\limits_{j=1}^{m}\,
\zeta_{_{\scriptstyle j}}(x) \wedge
\gamma_{_{\scriptstyle j}}(x) \Bigr) =
d\Bigl( d \alpha(x)_
{\displaystyle |_{\scriptstyle {\rm (ED)}}} +
\sum\limits_{j=1}^{m}\,
\zeta_{_{\scriptstyle j}}(x)\wedge
\gamma_{_{\scriptstyle j}}(x) \Bigr),
\hfill
$
\\[2ex]
we see that from the condition (2.5), we get the condition (2.1).

We prove the second assertion by contradiction.
Let the gaps $\Theta_{\tau}^{},\,\tau=1,\ldots, s,$ on the domain $G^{\;\!\prime}$
have the nonzero total index
$
\sum\limits_{\tau= 1}^{s}\,
{\rm ind}_{{}_{\scriptstyle d\theta}}
\Theta_{{}_{\scriptstyle \tau}} \ne 0
$
with respect to the $(n-1)\!$-form $d\theta.$
Then the integral
$
\int\limits_{\partial \Xi} d\theta(x) \ne 0
$
\vspace{0.35ex}
on any compact hypersurface $\partial \Xi$ from the domain $G^{\;\!\prime}$
surrounding all gaps $\Theta_{{}_{\scriptstyle \tau}},\,\tau = 1,\ldots, s,$ but no other gaps.
\vspace{0.25ex}
Suppose that one of the compact hypersurfaces
\vspace{0.35ex}
$\partial\Xi$ is an compact integral hypersurface of
the system of exterior differential equations (ED); we denote it by $\partial\Sigma.$
Then the integral
$
\int\limits_{\partial\Sigma} d\theta(x) \ne 0.
$

On the other hand, the $(n-2)\!$-form $\theta$ is invariant on the domain $G^{\;\!\prime}$
\vspace{0.35ex}
with respect to the system of exterior differential equations (ED). Therefore the integral
$
\int\limits_{\partial\Sigma}d\theta(x)=0.
$
The obtained contradiction completes the proof of Theorem 2.2. \k
\\[4.25ex]
\centerline{
{\bf 2.2. Pfaff system of equations}
}
\\[1.75ex]
\indent
{\bf 2.2.1.\! The first boundedness test for the number of compact integral hypersurfaces.}
The following Theorems 2.3 and 2.4 \,[9, 16] are straightforward corollaries of Theorems 2.1 and 2.2 for a Pfaff system of equations.
\vspace{0.35ex}

{\bf Theorem 2.3}
({\sl
the first boundedness test for the number of compact integral hypersurfaces of a Pfaff system of equations}).
\vspace{0.35ex}
{\it
Let a domain $G^{\;\!\prime}\subset G$ have the homotopy group $\pi_{n-1}^{}(G^{\;\!\prime})$ of rank $r,$ and
let there exist $(n-2)\!$-forms $\alpha\in C^{\,2}(G^{\;\!\prime})$ and $\ell_{j}^{}\in C^{\,1}(G^{\;\!\prime}),\, j= 1,\ldots,m,$
\vspace{0.35ex}
such that on the domain $G^{\;\!\prime}$ the exterior differential of sum
\\[2ex]
\mbox{}\hfill
$
\displaystyle
d\Bigl( d\alpha(x)_{\displaystyle |_{\scriptstyle {\rm (Pf)}}} +
\sum_{j=1}^{m}\,
\omega_{_{\scriptstyle j}}(x) \wedge
\ell_{_{\scriptstyle j}}(x)\Bigr) =
B(x)\, dx_1^{} \wedge\ldots \wedge dx_n^{}\;\!,
\hfill
$
\\[2.25ex]
where the function $B\colon G^{\;\!\prime}\to \R$ is constant sign.
\vspace{0.25ex}
Then the Pfaff system of equations {\rm (Pf)} with $\omega_{_{\scriptstyle j}}\in C^{\infty}(G),\, j= 1,\ldots,m,$
\vspace{0.75ex}
in the domain $G^{\;\!\prime}$ has at most $r$ compact integral hypersurfaces.
}

{\bf Theorem 2.4.}
\vspace{0.35ex}
{\it
Let a domain $G^{\;\!\prime}\subset G$ have the homotopy group $\pi_{n-1}(G^{\;\!\prime})$ of rank $r,$
and let there exist $(n-2)\!$-forms
\vspace{0.25ex}
$\alpha \in C^{\,2}(G^{\;\!\prime}),\ \theta\in C^{\,2}(G^{\;\!\prime}),
\ \ell_j^{}\in C^{\,1}(G^{\;\!\prime}), \ j= 1,\ldots,m,$
such that the $(n-2)\!$-form $\theta$ is invariant on the domain $G^{\;\!\prime}$
\vspace{0.25ex}
with respect to the Pfaff system of equations {\rm (Pf)} with $\omega_{_{\scriptstyle j}}\in C^{\,\infty}(G),\, j= 1,\ldots,m,$
and the exterior differential
\\[2ex]
\mbox{}\hfill
$
\displaystyle
d\Bigl( d\alpha(x)_{\displaystyle |_{\scriptstyle
{\rm (Pf)}}} + d\theta(x) +
\sum_{j=1}^{m}\,
\omega_{_{\scriptstyle j}}(x) \wedge
\ell_{_{\scriptstyle j}}(x) \Bigr) =
B(x)\, dx_1^{} \wedge \ldots \wedge dx_n^{}
$
\ for all $x\in G^{\;\!\prime},
\hfill
$
\\[2ex]
where the function $B\colon G^{\;\!\prime}\to \R$ is constant sign.
Then the following assertions are valid{\rm:}
\vspace{0.25ex}

{\rm1)} the system {\rm (Pf)} has at most $r$ compact integral hypersurfaces in the domain $G^{\;\!\prime};$
\vspace{0.25ex}

{\rm2)}
any set of gaps of the domain $G^{\;\!\prime}$ surrounded by
a compact integral hypersurface of system {\rm (Pf)} has the zero total index with respect to the $(n-1)\!$-form $d\theta.$
}
\vspace{0.75ex}

{\bf Example 2.2.}
Let us consider the Pfaff system of equations
\\[1.5ex]
\mbox{}\hfill                       
$
\omega_1^{}(x)=0,
\quad
\omega_2^{}(x)=0,
$
\hfill (2.6)
\\[1ex]
where the differential 1-forms
\\[1.5ex]
\mbox{}\hfill
$
\omega_1^{}(x)=
x_1^{}\,dx_1^{} + x_2^{}\,dx_2^{} + g(x)(x_4^{}\,dx_3^{}- x_3^{}\,dx_4^{})
$
\ for all $x\in\R^4,
\hfill
$
\\[2ex]
\mbox{}\hfill
$
\omega_2^{}(x)=
x_1^{}\,dx_1^{} + x_2^{}\,dx_2^{} + (2x_3^{} - x_4^{})\,dx_3^{} + (x_3^{} + 2x_4^{})\,dx_4^{}
$
\ for all
$
x\in\R^4,
\hfill
$
\\[2ex]
the scalar function
\vspace{0.5ex}
$
g(x)=x_1^2 + x_2^2 + x_3^2 + x_4^2
$
for all
$
x\in\R^4.
$

The Pfaff system (2.6) such that the differential 2-forms
$
\alpha(x) = 0
$
for all
$
x \in \R^4,
$
\\[1.75ex]
\mbox{}\hfill
$
\ell_{1}^{}(x) = g^{{}-1}(x)\,dx_1^{} \wedge dx_2^{}
$
\ for all $x \in G^{\;\!\prime},
\qquad
\ell_2^{} (x) = dx_3^{} \wedge dx_4^{}
$
\ for all
$
x\in\R^4
\hfill
$
\\[2ex]
on the domain $G^{\;\!\prime}=\R^4\backslash\{(0,0,0,0)\}$ satisfy the relations
\\[2ex]
\mbox{}\hfill
$
\omega_{1}^{}(x) \wedge \ell_{1}^{}(x) =
x_4^{}\,dx_1^{} \wedge dx_2^{} \wedge dx_3^{} - x_3^{}\,dx_1^{} \wedge dx_2^{} \wedge dx_4^{}\;\!,
\hfill
$
\\[2ex]
\mbox{}\hfill
$
\omega_{2}^{}(x) \wedge \ell_{2}^{}(x) =
x_1^{}\,dx_1^{} \wedge dx_3^{} \wedge dx_4^{} + x_2^{}\,dx_2^{} \wedge dx_3^{} \wedge dx_4^{}\;\!,
\hfill
$
\\[2.25ex]
\mbox{}\hfill
$
d\bigl(
\omega_{1}^{}(x) \wedge \ell_{1}^{}(x) +
\omega_{2}^{}(x) \wedge \ell_{2}^{}(x)
\bigr) =
{}-2\,dx_1^{} \wedge dx_2^{} \wedge dx_3^{} \wedge dx_4^{}.
\hfill
$
\\[2.15ex]
\indent
Therefore (by Theorem 2.3) in the domain $G^{\;\!\prime}=\R^4\backslash\{(0,0,0,0)\}$
\vspace{0.25ex}
with homotopy group $\pi_{3}^{}(G^{\;\!\prime})$  of rank 1 the Pfaff system (2.6) can have at most one compact integral hypersurface.
\vspace{0.25ex}

Now, if we take into account the fact that
\\[1.5ex]
\mbox{}\hfill
$
d(g(x)-1) = \bigl(\omega_1^{}(x) +\omega_2^{}(x)\bigr)_{
\displaystyle |_{\scriptstyle g(x) = 1}}\;\!,
\hfill
$
\\[1.75ex]
then we see that the sphere
\vspace{0.25ex}
$
S^3 = \{x\colon g(x) = 1\}
$
is this single compact integral hypersurface of the Pfaff system of equations (2.6).
\vspace{1ex}

{\bf Example 2.3.}
Consider the Pfaff system of equations
\\[1.5ex]
\mbox{}\hfill                    
$
\omega_1^{}(x)=0,
\quad
\omega_2^{}(x)=0,
$
\hfill (2.7)
\\[1.5ex]
where the differential 1-forms
\\[1.5ex]
\mbox{}\hfill
$
\omega_1^{}(x)=
x_3^{}\bigl((x_1^{} - 2)^2 + x_2^2 + x_3^2\bigr)\!\;dx_1^{}+ x_3^{}\,dx_2^{} +
x_2^{}\bigl((x_1^{} -2)^2 + x_2^2 + x_3^2\bigr)\!\;dx_3^{}
$
\ for all
$
x\in \R^{3},
\hfill
$
\\[2.25ex]
\mbox{}\hfill
$
\omega_2^{}(x)=
\Bigl(\bigl(\sqrt{x_1^2 + x_2^2} - 2\bigr)^2 + x_3^2 - 1 +
x_1^{}\sqrt{x_1^2 + x_2^2}\,\bigl(\sqrt{x_1^2 + x_2^2} -2\bigr)\Bigr)\;\!dx_1^{} \ +
\hfill
$
\\[2ex]
\mbox{}\hfill
$
+ \Bigl(\bigl(\sqrt{x_1^2 + x_2^2} - 2\bigr)^2 + x_3^2 - 1 + x_2^{}
\sqrt{x_1^2 + x_2^2}\,
\bigl(\sqrt{x_1^2 + x_2^2} - 2\bigr)\Bigr)\;\!dx_2^{} \ +
\hfill
$
\\[2ex]
\mbox{}\hfill
$
+ \ x_3^{}(x_1^2 + x_2^2)\,dx_3^{}
$
\ for all
$
x\in \R^{3}.
\hfill
$
\\[2.25ex]
\indent
The Pfaff system (2.7) such that the differential 1-forms
$
\alpha(x) = 0
$
for all
$x\in\! \R^{3},
$
\\[2.25ex]
\mbox{}\hfill
$
\ell_{1}^{}(x) =\dfrac{1}{(x_1^{} - 2)^2 + x_2^2 + x_3^2}\ dx_2^{}
$
\ for all $x\in G^{\;\!\prime},
\qquad
\ell_{2}^{}(x) = 0
$
\ for all
$x\in \R^{3},
\hfill
$
\\[2.25ex]
on the domain $G^{\;\!\prime}=\R^3\backslash\{(2,0,0)\}$ satisfy the relations
\\[1.75ex]
\mbox{}\hfill
$
\omega_{1}^{}(x) \wedge \ell_{1}^{}(x) =
x_3^{}\,dx_1^{} \wedge dx_2^{} - x_2^{}\, dx_2^{}\, \wedge dx_3^{},
\qquad
\omega_{2}^{}(x)\wedge \ell_{2}^{}(x) = 0,
\hfill
$
\\[2.25ex]
\mbox{}\hfill
$
d
\bigl(\omega_{1}^{}(x)\wedge \ell_{1}^{}(x) +
\omega_{2}^{}(x)\wedge \ell_{2}^{}(x)\bigr)
= dx_1^{}\wedge dx_2^{} \wedge dx_3^{}.
\hfill
$
\\[2ex]
\indent
Consequently (by Theorem 2.3) in the domain $G^{\;\!\prime}=\R^3\backslash\{(2,0,0)\}$
\vspace{0.25ex}
with homotopy group $\pi_{2}(G^{\;\!\prime})$  of rank 1 the Pfaff system (2.7) can have at most one compact integral hypersurface.
\vspace{0.35ex}

Since
\\[1ex]
\mbox{}\hfill
$
dw(x) = \dfrac{2}{x_1^2 + x_2^2}\
\omega_{2}^{}(x)_{\displaystyle |_{\scriptstyle w(x)= 0}}\,,
\hfill
$
\\[1.75ex]
where the scalar function
\\[1.75ex]
\mbox{}\hfill
$
w\colon x\to\bigl(\sqrt{x_1^2 + x_2^2} - 2\bigr)^2 + x_3^2 - 1
$
\ for all
$x \in \R^3,
\hfill
$
\\[2ex]
we see that the two-dimensional torus
$
\{x\colon w(x) = 0\}
$
\vspace{0.35ex}
is this single compact integral hypersurface of the Pfaff system of equations (2.7).
\vspace{1ex}

{\bf Theorem 2.5.}
{\it
A linear Pfaff system has no isolated compact integral hypersurfaces.
}
\vspace{0.5ex}

{\sl Proof.}
Suppose a linear Pfaff system of equations is induced the 1-forms
\vspace{0.35ex}
$\omega_j^{},\ j=1,\ldots,m,$ with linear coordinate functions.
Let us consider the following two logical possibilities:
\vspace{0.25ex}

1) there exists an index $k \in \{1,\ldots,m\}$ such that
$
d\omega_k^{}(x) \ne 0;
$

2) $d\omega_j^{}(x) = 0$ for all $x\in \R^n,\ j=1,\ldots,m.$
\vspace{0.5ex}

In the first case, the exterior differential has the form
\\[2ex]
\mbox{}\hfill
$
\displaystyle
d\omega_k^{} (x) =
\sum\limits_{1 \leq i < \tau \leq n}
c_{i \tau k}^{}\,
dx_i^{} \wedge dx_{\tau}^{},
\hfill
$
\\[2ex]
where the coefficients $c_{i \tau k}^{}$ are real numbers that do not vanish simultaneously.
Let $c_{{}_{\scriptstyle \lambda \rho k}} \ne 0,$ $\lambda < \rho.$
We take on the space $\R^n$ the $(n-2)\!$-form
\\[2ex]
\mbox{}\hfill
$
\ell(x) = dx_1^{}\wedge \ldots \wedge
dx_{{}_{\scriptstyle \lambda -1}}  \wedge
dx_{{}_{\scriptstyle \lambda +1}} \wedge
dx_{{}_{\scriptstyle \lambda +2}} \wedge \ldots  \wedge
dx_{{}_{\scriptstyle \rho -1}}\, \wedge
dx_{{}_{\scriptstyle \rho + 1}} \wedge
dx_{{}_{\scriptstyle \rho + 2}} \wedge \ldots  \wedge
dx_{{}_{\scriptstyle n}}\,.
\hfill
$
\\[2ex]
\indent
Then the exterior differential of exterior product
\\[2ex]
\mbox{}\hfill
$
d\bigl(
\omega_k^{} (x) \wedge \ell(x) \bigr) =
{}\pm c_{{}_{\scriptstyle \lambda \rho k}}\,
dx_1^{} \wedge \ldots \wedge dx_n^{}
$
\ for all
$
x \in \R^n.
\hfill
$
\\[2ex]
\indent
Using Theorem 2.3, we obtain the linear Pfaff system of equations has neither isolated compact integral hypersurfaces
nor compact integral hypersurfaces.

In the second case, the 1-forms $\omega_{j}^{},\, j=1,\ldots,m,$
\vspace{0.25ex}
are total differentials on the space $\R^n.$
Therefore the linear Pfaff system has a basis of first integrals
\vspace{0.35ex}
$
F_j^{}\colon \R^n \to \R,\  j=1,\ldots,m,
$
where
$
F_j^{}
$
\vspace{0.25ex}
are some polynomials of degree
$
{\rm deg}\;\! F_j^{}\leq 2,\  j=1,\ldots,m.
$
Since this is an algebraic basis, it follows that there is no isolated compact integral hypersurface.  \k
\vspace{0.5ex}

From the proof of Theorem 2.5, we get the following corollaries about
compact integral hypersurfaces of a linear Pfaff system of equations.
\vspace{0.75ex}

{\bf Corollary 2.1.}
{\it
If for the linear Pfaff system of equations
\vspace{0.35ex}
$
\omega_j^{}(x) = 0, \
j=1,\ldots,m,
$
there exists an index $k \in \{1,\ldots,m\}$ such that
$
d\omega_k^{}(x) \ne 0,
$
\vspace{0.25ex}
then this linear Pfaff system
has no compact integral hypersurfaces.
}
\vspace{0.75ex}

{\bf Corollary 2.2.}
{\it
An non-isolated compact integral hypersurfaces
of a linear Pfaff system of equations are
an algebraic hypersurfaces of second order.
}
\\[3.25ex]
\indent
{\bf
2.2.2.\!
The second boundedness test for the number of compact integral hypersurfaces.}
\\[0.5ex]
\indent
{\bf Theorem 2.6}
({\sl the second boundedness test for the number of compact integral hypersurfaces of a Pfaff system of equations}).
\vspace{0.25ex}
{\it
\!\!Let a domain $\!G^{\;\!\prime}\!\subset\! G\!$  have the homotopy group $\!\pi_{n-1}^{}(G^{\;\!\prime})$ of rank $r,$
let there exists a vector field $V\in C^1(G^{\;\!\prime})$
orthogonal to the vector fields
\\[1.75ex]
\mbox{}\hfill
$
W_j^{}\colon x\to \bigl(w_{_{\scriptstyle j1}}(x),\ldots,w_{_{\scriptstyle jn}}(x)\bigr)
$
\ for all
$
x\in G^{\;\!\prime},
\quad
j=1,\ldots,m,
\hfill
$
\\[1.75ex]
such that
\vspace{0.25ex}
${\rm div}\,V$ is constant sign on the domain $G^{\;\!\prime}.$
Then the Pfaff system {\rm (Pf)} with $\omega_j^{}\!\in\! C^{\infty}(G),\, j\!=\!1,\ldots,m,$
\vspace{0.75ex}
in the domain  $G^{\;\!\prime}$ has at most $r$ compact integral hypersurfaces.
}

The {\sl proof} of this theorem is similar to that of  Theorem 1.1 when $\nu = n$ [9, 16]
and is based on the following
\vspace{0.5ex}

{\bf Lemma 2.2.}
{\it
Let the assumptions of Theorem {\rm 2.6} be valid.
Then in a subdomain $\Omega$ of the domain $G^{\;\!\prime}\subset G$
with homotopy group $\pi_{n-1}^{}(\Omega)$ of rank $s,\ s \leq r,$
\vspace{0.25ex}
the situation described in Lemma {\rm 2.1} is impossible
\vspace{0.25ex}
for compact integral hypersurfaces of the Pfaff system of equations {\rm (Pf)}
with $\omega_j\in C^{\infty}(G),\ j=1,\ldots,m.$
}
\vspace{0.5ex}

{\sl Proof.}
The proof of Lemma 2.2 is similar to that of Lemma 2.1, but in this case we use the Ostrogradskii formula
\\[2ex]
\mbox{}\hfill
$
\displaystyle
\int\limits_{\partial \Sigma} V(x) n(x)\, dS =
({}-1)^n \int\limits_\Sigma {\rm div}\, V(x)\, d\Sigma,
\hfill
$
\\[1.5ex]
where $n$ is the unit normal on the integral hypersurface
$
\partial \Sigma = \bigcup\limits_{\tau=1}^{s+1}
\partial\Sigma_{\tau}
$
of system (Pf), and take into account the orthogonality of the vector fields
\vspace{0.5ex}
$W_j^{},\ j=1,\ldots,m,$
to the hypersurface $\partial \Sigma$ and the constant sign of the scalar function
${\rm div}\,V$ on the domain $G^{\;\!\prime}.$ \k
\vspace{0.5ex}

{\bf Example 2.4.}
Consider the Pfaff system of equations
\\[1.5ex]
\mbox{}\hfill                    
$
\omega_1^{}(x)=0,
\qquad
\omega_2^{}(x)=0,
$
\hfill (2.8)
\\[1.5ex]
where the differential 1-forms
\\[1.5ex]
\mbox{}\hfill
$
\omega_1^{}(x)=\bigl(x_1^{} - x_2^{} + x_2^{}\;\!g(x)\bigr)\;\!dx_1^{} +
\bigl(x_1^{} + x_2^{} - x_1^{}\;\!g(x)\bigr)\;\!dx_2^{} \ +
\hfill
$
\\[1.75ex]
\mbox{}\hfill
$
+\
\bigl(x_3^{} - x_4^{} + x_4^{}\;\!g(x)\bigr)\;\!dx_3^{} +
\bigl(x_3^{} + x_4^{} - x_3^{}\;\!g(x)\bigr)\;\!dx_4^{}
$
\, \ for all
$
x\in\R^4,
\hfill
$
\\[2ex]
\mbox{}\hfill
$
\omega_2^{}(x)=
\bigl(x_3^{} - x_4^{} + x_4^{}\;\!g(x)\bigr)\;\!dx_1^{} +
\bigl(x_3^{} + x_4^{} - x_3^{}\;\!g(x)\bigr)\;\!dx_2^{} \ +
\hfill
$
\\[1.75ex]
\mbox{}\hfill
$
+\
\bigl(x_1^{} - x_2^{} + x_2^{}\;\!g(x)\bigr)\;\!dx_3^{} +
\bigl(x_1^{} + x_2^{} - x_1^{}\;\!g(x)\bigr)\;\!dx_4^{}
$
\, \ for all
$
x\in\R^4,
\hfill
$
\\[2ex]
the scalar function
$
g\colon x\to x_1^2 + x_2^2 + x_3^2 + x_4^2
$
for all $x\in\R^4.$
\vspace{0.75ex}

The continuously differentiable on the domain $G^{\;\!\prime}=\R^4\backslash\{(0,0,0,0)\}$ vector field
\\[1.5ex]
\mbox{}\hfill
$
V\colon x \to g^{{}-3}(x)
\bigl({}-x_1 - x_2 + x_1g(x),\, x_1 - x_2 + x_2g(x),
{}-x_3 - x_4 + x_3g(x),\, x_3 - x_4 + x_4g(x)\bigr)
\hfill
$
\\[1.5ex]
is orthogonal to the vector fields $W_1^{}$ and $W_2^{}$ induced by the system (2.8).
\vspace{0.35ex}
Moreover, the function
$
{\rm div}\,V(x) = 2g^{{}-3}(x)$ for all $x\in G^{\;\!\prime}
$
is positive on the domain $G^{\;\!\prime}.$
\vspace{0.5ex}

By Theorem 2.6,
\vspace{0.35ex}
the Pfaff system (2.8) in the domain $G^{\;\!\prime}=\R^4\backslash\{(0,0,0,0)\}$ with
the homotopy group $\pi_{3}(G^{\;\!\prime})$ of rank 1 can have at most one compact integral hypersurface.
\vspace{0.65ex}

Taking into account
\vspace{0.65ex}
$
d(g(x)-1) = 2\omega_1(x)_{\displaystyle |_{\scriptstyle  g(x) = 1}},
$
we obtain the sphere
$
S^3 = \{x\colon g(x) = 1\}
$
is this single compact integral hypersurface of the Pfaff system of equations (2.8).
\vspace{1ex}

Let us consider a Pfaffian equation
\\[1.5ex]
\mbox{}\hfill   
$
\omega(x) = 0,
$
\hfill (2.9)
\\[1.25ex]
where the differential 1-form
\\[1.75ex]
\mbox{}\hfill
$
\displaystyle
\omega(x)  = \sum\limits_{i=1}^{n} a_i^{}(x)\, dx_i^{}
$
\ for all $x\in G,
\hfill
$
\\[1.5ex]
the scalar functions $a_i^{}\in C^{\infty}(G), \ i = 1,\ldots, n,$
the domain $G\subset \R^n.$
\vspace{0.35ex}

We assume that the Pfaffian equation (2.9)
\vspace{0.25ex}
is completely integrable [17, 18, 19] on the do\-ma\-in $G,$
i.e., the Frobenius condition
$d\;\!\omega(x) \wedge\, \omega(x) = 0$ for all $x\in G$ is satisfied.
\vspace{0.35ex}

The following corollary [16] is an immediate consequence of Theorem 2.6.
\vspace{0.5ex}

{\bf Corollary 2.3.}
{\it
Let a domain $G^{\;\!\prime}\subset G$
\vspace{0.25ex}
have the homotopy group $\pi_{n-1}^{}(G^{\;\!\prime})$ of rank $ r,$
let there exists a vector field $V\in C^1(G^{\;\!\prime})$ orthogonal on the domain $G^{\;\!\prime}$
to the vector field
\\[1.5ex]
\mbox{}\hfill
$
A\colon x\to \bigl(a_1^{}(x),\ldots,a_n^{}(x)\bigr)
$
\ for all $x\in G
\hfill
$
\\[1.5ex]
such that the divergence ${\rm div}\,V$  is constant sign on  the domain $G^{\;\!\prime}.$
Then the completely integrable Pfaffian equation {\rm (2.9)} has at most $r$ compact leaves in the domain $G^{\;\!\prime}.$
}
\vspace{0.5ex}

For the practical use of Corollary 2.3 we construct the set of $n-1$ linearly independent vector fields
\\[1.5ex]
\mbox{}\hfill   
$
{\Phi}_{\tau}^{}\colon x\to\,
 \left(\varphi_1^{}(x), \ldots, \varphi_n^{}(x)\right)
$
\ for all $x \in G,
\quad \
\varphi_{l}^{}\colon x\to\,  \delta_{l \tau }^{}\;\! a_n^{}(x)
$
\ for all $x \in G,
\hfill
$
\\[0.5ex]
\mbox{}\hfill (2.10)
\\
\mbox{}\hfill
$
\varphi_n^{}\colon x\to  {}- a_{\tau}^{}(x)$
\ for all $x \in G,
\quad
l=1,\ldots,n-1, \
\tau=1,\ldots,n-1,
\hfill
$
\\[2ex]
where $\delta_{l\;\!\tau}^{}$ is the Kronecker delta.
\vspace{0.35ex}
The vector fields (2.10) are orthogonal to the vector field
$A\colon x\to (a_1^{}(x),\ldots,a_n^{}(x))$ for all $x \in G^{\;\!\prime}.$
\vspace{0.35ex}
Therefore any vector field $V\colon G^{\;\!\prime}\to \R^n$
orthogonal to $A$ on the domain $G^{\;\!\prime}$ can be represented in the basis
${\Phi}_{\tau}^{},\  \tau = 1,\ldots,n-1,$ as the sum
\\[1.5ex]
\mbox{}\hfill
$
\displaystyle
V(x) = \sum\limits_{\tau=1}^{n-1} g_{\tau}^{}(x)\;\!{\Phi}_{\tau}^{}(x)
$
\ for all
$x \in G^{\;\!\prime},
\hfill
$
\\[1.5ex]
where the scalar functions $g_{\tau}^{}\in C^1(G^{\;\!\prime}),\  \tau = 1,\ldots,n-1.$
Now using Corollary 2.3, we get
\vspace{0.75ex}

{\bf Theorem 2.7.}
{\it
Let a domain
\vspace{0.35ex}
$G^{\;\!\prime}\subset G$ have the homotopy group $\pi_{n-1}^{}(G^{\;\!\prime})$ of rank $r,$
and let there exist scalar functions $g_{\tau}^{}\in C^1(G^{\;\!\prime}), \ \tau = 1,\ldots,n-1,$
such that the vector field
\\[1.5ex]
\mbox{}\hfill
$
\displaystyle
V(x) = \sum\limits_{\tau=1}^{n-1} g_{\tau}^{}(x){\Phi}_{\tau}^{}(x)
$
\ for all
$x \in G^{\;\!\prime}
\hfill
$
\\[1.5ex]
{\rm(}which is constructed on the basis of the linearly independent vector fields {\rm(2.10))}
has the constant sign divergence ${\rm div}\;\!V$ on the domain $G^{\;\!\prime}.$
Then the co\-m\-p\-le\-te\-ly integrable Pfaffian equation {\rm (2.9)} has at most $r$ compact leaves in the domain $G^{\;\!\prime}.$
}
\vspace{0.75ex}

Let us indicate tests for the absence [16] of compact leaves of the Pfaffian equation (2.9).
\vspace{0.75ex}

{\bf  Theorem 2.8.}
\vspace{0.25ex}
{\it
Suppose a domain $G^{\;\!\prime}\subset G$ has the  trivial homotopy group $\pi_{n-1}^{}(G^{\;\!\prime}),$
there exists a scalar function $\mu\in C^1(G^{\;\!\prime})$ such that the vector field
\\[1.5ex]
\mbox{}\hfill
$
\mu\;\!A\colon x\to \mu(x)(a_1^{}(x),\ldots,a_n^{}(x))
$
\ for all $x\in G^{\;\!\prime}
\hfill
$
\\[1.5ex]
can vanish only on a null set of the $(n-2)\!$-dimensional measure,
and the vector field $\mu\;\!A$ is solenoidal on the domain $G^{\;\!\prime},$ i.e.,
\\[1.5ex]
\mbox{}\hfill
$
{\rm div} \left(\mu(x)\;\!A(x)\right)  =  0$
\ for all
$
x \in G^{\;\!\prime}.
\hfill
$
\\[1.5ex]
Then the completely integrable\!
\vspace{0.5ex}
Pfaffian\! equation\! {\rm (2.9)}\! has no compact leaves in the domain $\!G^{\;\!\prime}\!.$
}

{\sl Proof.}
The proof is performed similarly to that of Theorem 1.1 when $m=1,\, \nu =n$ and is based on the Ostrogradskii formula
\\[2ex]
\mbox{}\hfill
$
\displaystyle
\int\limits_L
(\mu(x)\;\!A(x),  n(x))\,dS =
({}-1)^n
\int\limits_{\Xi} \mbox{\rm div}\;\!\left(\mu(x)\;\!{A}(x)\right)\, d\Xi,
\hfill
$
\\[1.5ex]
where $(\cdot,\cdot)$ is the operation of scalar product of vectors,
$n$ is the unit external normal vector to the leaf $L=\partial \Xi.$
In our case this formula is fails.

Namely, since the field $\mu\,{A}$ is orthogonal to the leaf  $L,$ it follows that the integral
\\[1.5ex]
\mbox{}\hfill
$
\displaystyle
\int\limits_L (\mu(x)\;\! {A}(x),  n(x))\,dS\ne 0.
\hfill
$
\\[1.5ex]
But on the other hand, since the vector field $\mu A$ is solenoidal, we have
\\[1.5ex]
\mbox{}\hfill
$
\displaystyle
\int\limits_{\Xi} \mbox{\rm div} \; \left(\mu(x){A}(x)\right)\; d\Xi=0. \ \k
\hfill
$
\\[2ex]
\indent
{\bf  Corollary 2.4.}
{\it
If the vector field
\\[1.5ex]
\mbox{}\hfill
$
A\colon x\to (a_1^{}(x),\ldots,a_n^{}(x) )
$
\ for all
$x\in G^{\;\!\prime}
\hfill
$
\\[1.5ex]
is solenoidal on the domain $G^{\;\!\prime}\subset G $ with trivial homotopy group
\vspace{0.35ex}
$\pi_{n-1}^{}(G^{\;\!\prime})$ and can vanish only on a null set of the $(n-2)\!$-dimensional measure,
\vspace{0.25ex}
then the completely integrable Pfaffian equation {\rm (2.9)} has no compact leaves in the domain $G^{\;\!\prime}.$
}
\vspace{0.75ex}

{\bf Example 2.5.}
The Pfaffian equation
\\[1.75ex]
\mbox{}\hfill                    
$
yz\, dx+ 2xz\, dy+ 3xy\, dz=0
$
\hfill (2.11)
\\[1.75ex]
is completely integrable on $\R^3.$ Indeed,
the vector field
\\[1.5ex]
\mbox{}\hfill
$
A\colon (x,y,z)\to\, (yz, 2xz, 3xy)
$
\ for all
$(x,y,z)\in\R^3
\hfill
$
\\[1.5ex]
and its rotor
\\[0.5ex]
\mbox{}\hfill
$
{\rm rot}\,A(x,y,z) =(x,{}- 2y, z)$ for all $(x,y,z)\in\R^3
\hfill
$
\\[1.75ex]
are orthogonal, i.e., $(A(x,y,z), {\rm rot}\,A(x,y,z) ) =0$ for all $(x,y,z)\in\R^3.$
\vspace{0.5ex}

The vector field $A$ is  solenoidal, i.e.,
${\rm div}\,A(x,y,z) =0$ for all $(x,y,z)\in\R^3.$
\vspace{0.5ex}

By Corollary 2.4, the completely integrable Pfaffian equation (2.11) has no compact leaves.

Note that
\\[0.5ex]
\mbox{}\hfill
$
d(xy^2z^3)_{|_{(2.11)}}=0
$
\ for all
$(x,y,z)\in\R^3.
\hfill
$
\\[1.5ex]
Therefore the Pfaffian equation (2.11) is determined the foliation $xy^2z^3=C$ of space $\R^3.$
This foliation doesn't contain compact leaves.
\\[4.25ex]
\centerline{
{\bf 2.3. Ordinary autonomous differential system}
}
\\[1.75ex]
\indent
{\bf 2.3.1.\! The first boundedness test for the number of compact integral hypersurfaces.}
The  ordinary autonomous differential system {\rm (D)} is induced the
system of  $n(n - 1)/2$ Pfaffian equations
\\[1.5ex]
\mbox{}\hfill                                         
$
\psi_{{}_{\scriptstyle qh}}(x) = 0,
\quad
1 \leq q < h \leq n,
$
\hfill (2.12)
\\[1.5ex]
where the 1-forms
\\[1.5ex]
\mbox{}\hfill
$
\psi_{{}_{\scriptstyle qh}}(x) =
f_{{}_{\scriptstyle q}}(x)\,dx_h^{} -
f_h^{}(x)\,dx_{{}_{\scriptstyle q}}
$
\ for all $
x\in G,
\quad
1\leq q < h \leq n,
\hfill
$
\\[1.75ex]
are closed on the domain $G.$

An autonomous integral basis of system {\rm (D)} is a basis of first integrals of the Pfaff system of equations (2.12), and vice versa.
This allows one to generalize Theorems 2.3 and 2.4 to the case of system {\rm (D)}, and we call
Theorems 2.3{\rm D}
{\sl {\rm(}the first boundedness test for the number of compact integral hypersurfaces of an ordinary autonomous differential system})
and 2.4{\rm D}.

As to Theorem 2.4{\rm D} in the case of system (0.1), we note that a scalar function $\theta\colon G\to \R$ invariant with respect
to the Pfaff system of equations (2.12) induced by system (0.1) is an autonomous general integral of system (0.1).
\vspace{0.5ex}

{\bf Example 2.6.}
Using the first boundedness test for the number of compact integral hy\-per\-sur\-fa\-ces (Theorem 2.3{\rm D}),
we could prove that the ordinary autonomous differential sys\-tem
\\[2.5ex]
\mbox{}\hfill                                            
$
\dfrac{dx_1^{}}{dt} = x_3^{}\;\!g(x),
\quad \
\dfrac{dx_2^{}}{dt} = {}-x_2^{} + x_3^{} + x_2^{}\;\!g(x),
\quad \
\dfrac{dx_3^{}}{dt} = {}-x_2^{} - x_1^{}\;\!g(x),
$
\hfill (2.13)
\\[2.5ex]
where the scalar function
\vspace{0.35ex}
$g\colon x\to x_1^2 +x_2^2 + x_3^2$ for all $x\in\R^3,
$
has one compact integral surface in the domain $G^{\;\!\prime} =\R^3\backslash\{(0,0,0)\}.$
\vspace{0.5ex}

We claim that the sphere
$
S^2 = \{x\colon g(x) = 1\}
$
is a compact integral surface of the dif\-fe\-ren\-ti\-al system (2.13). Indeed,
\\[2ex]
\mbox{}\hfill
$
\dfrac{d ( g(x)-1 )}{dt}_{\bigl|_{\scriptstyle (2.13)}} =  2x_2^2\,(g(x)-1)$
\ for all $x \in \R^3.
\hfill
$
\\[2ex]
\indent
Using the ordinary differential sys\-tem (2.13), we form the Pfaffian equations
\\[1.5ex]
\mbox{}\hfill
$
\psi_{{}_{\scriptstyle 12}}(x) = 0,
\hfill
$
\\[1.5ex]
where
$
\psi_{{}_{\scriptstyle 12}}(x) =
x_3^{}\;\!g(x)\,dx_2^{} - \bigl({}-x_2^{} + x_3^{} + x_2^{}\;\!g(x)\bigr)\;\!dx_1^{}$
for all $x \in \R^3.$
\vspace{0.75ex}

The continuously differentiable on the domain  $G^{\;\!\prime}$ 1-form
\\[1.5ex]
\mbox{}\hfill
$
\omega(x) = {}-\dfrac{1}{g(x)}\ dx_1^{}
$
\ for all
$
x \in G^{\;\!\prime}
\hfill
$
\\[1.5ex]
such that the exterior differential of exterior product
\\[1.5ex]
\mbox{}\hfill
$
d\bigl(\psi_{{}_{\scriptstyle 12}}(x) \wedge\omega(x)\bigr) =
dx_1^{} \wedge dx_2^{}\wedge dx_{3}^{}
$
\ for all
$
x \in G^{\;\!\prime}.
\hfill
$
\\[1.5ex]
\indent
Therefore, by Theorem 2.3{\rm D}, the system (2.13) in the domain
$G^{\;\!\prime}$ with homotopy group $\pi_{2}(G^{\;\!\prime})$ of ranf 1 has at most one compact integral surface.
This surface is the sphere $S^2.$
\\[2.5ex]
\indent
{\bf 2.3.2.\! The second boundedness test for the number of compact integral hy\-per\-sur\-fa\-ces.}
Using the first boundedness test for the number of compact integral hy\-per\-sur\-fa\-ces
of an ordinary autonomous differential system (Theorem 2.3{\rm D}), we could prove Theorem 2.9.
In this assertion the boundedness for the number of compact integral hy\-per\-sur\-fa\-ces
is established by the form of systems {\rm (D)}
\vspace{0.5ex}
and is not used the Pfaff system of equations (2.12).

{\bf Theorem 2.9.}
{\it
Let a domain
\vspace{0.25ex}
$G^{\;\!\prime}\subset G$ have the homotopy group $\pi_{n-1}(G^{\;\!\prime})$ of rank $r,$
let there exists a scalar function $\varphi\in C^1(G^{\;\!\prime})$
\vspace{0.25ex}
such that the divergence ${\rm div}\,Z$ of the vector field
$\!
Z\colon x \to \varphi(x)f(x)\!$
for all
$\!
x\!\in\! G^{\;\!\prime}\!
$
\vspace{0.35ex}
is constant sign on the domain $\!G^{\;\!\prime}\!.$
Then the system {\rm (D)} with $f\in C^{\infty}(G)$  has at most $r$ compact integral hypersurfaces
\vspace{0.5ex}
in the domain $G^{\;\!\prime}.$
}

{\sl Proof.}
From Theorem 2.3{\rm D}, it follows that
if as the 1-forms $\omega_{i}^{},\ i = 1,\ldots, n,$ we choose
\\[1.75ex]
\mbox{}\hfill
$
\omega_{{}_{\scriptstyle \nu}}(x) =
f_{{}_{\scriptstyle \nu}}(x)\,dx_{{}_{\scriptstyle \nu + 1}} -
f_{{}_{\scriptstyle \nu + 1}}(x)\,dx_{{}_{\scriptstyle \nu}}
$
\ for all
$
x\in G,
\quad
\nu = 1,\ldots, n-1,
\hfill
$
\\[2.75ex]
\mbox{}\hfill
$
\omega_{{}_{\scriptstyle n}}(x) =
f_{{}_{\scriptstyle n}}(x)\,dx_{_{\scriptstyle 1}} -
f_{_{\scriptstyle 1}}(x)\,dx_{{}_{\scriptstyle n}}
$
\ for all
$
x\in G,
\hfill
$
\\[2ex]
and as the $(n - 2)\!$-forms $\ell_{i}^{}, \ i = 1,\ldots,n,$ we take
\\[1.75ex]
\mbox{}\hfill
$
\ell_{{}_{\scriptstyle \nu}}(x) =
\dfrac12\  \varphi(x)\,
dx_{{}_{\scriptstyle 1}}\wedge \ldots \wedge
dx_{{}_{\scriptstyle \nu -1}}\wedge
dx_{{}_{\scriptstyle \nu + 2}}\wedge
dx_{{}_{\scriptstyle \nu + 3}} \wedge\ldots\wedge
dx_{{}_{\scriptstyle n}}
$
\ for all
$
x\in G^{\;\!\prime},
\hfill
$
\\[2.75ex]
\mbox{}\hfill
$
\ell_{{}_{\scriptstyle n}}(x) =
(-1)^{n+1}\,\dfrac12\ \varphi(x)\,
dx_{{}_{\scriptstyle 2}}\wedge \ldots \wedge
dx_{{}_{\scriptstyle n-1}}
$
\ for all
$
x\in G^{\;\!\prime},
\quad
\nu = 1,\ldots, n-1,
\hfill
$
\\[1.75ex]
where the scalar function $\varphi\in C^1(G^{\;\!\prime}),$
then the sum of exterior products
\\[2ex]
\mbox{}\hfill
$
\displaystyle
\sum\limits_{i=1}^{n}\,
\omega_{i}^{}(x)\wedge \ell_{i}^{}(x) =
\hfill
$
\\[1.25ex]
\mbox{}\hfill
$
\displaystyle
=\sum\limits_{i=1}^{n}\, ({}-1)^{i + 1}
f_{{}_{\scriptstyle i}}(x)\,
dx_{{}_{\scriptstyle 1}}\wedge \ldots \wedge
dx_{{}_{\scriptstyle i-1}}\wedge
dx_{{}_{\scriptstyle i+1}}\wedge
dx_{{}_{\scriptstyle i+2}} \wedge \ldots\wedge
dx_{{}_{\scriptstyle n}}
$
\ for all $x\in G^{\;\!\prime}
\hfill
$
\\[2ex]
and the exterior differential
\\[2ex]
\mbox{}\hfill
$
\displaystyle
d\Bigl(\varphi(x)\sum\limits_{i=1}^{n}\,
\omega_i^{}(x) \wedge \ell_i^{}(x)\Bigr) =
{\rm div}\,Z(x)\,dx_{1}^{}\wedge \ldots
\wedge dx_{n}^{}
$
\ for all
$
x\in G^{\;\!\prime}.
\hfill
$
\\[2ex]
\indent
Consequently, by Theorem 2.3{\rm D}, Theorem 2.9 is valid. \k
\vspace{1ex}

Note also that Theorem 2.9 is a consequence of Theorem 2.6 on the case of
\vspace{0.35ex}
the ordinary autonomous differential system {\rm (D)} with $f\in C^{\infty}(G).$
It follows from that
the vector field
\\[1.75ex]
\mbox{}\hfill
$
V\colon x \to\  \varphi(x)\;\!f(x)
$
\ for all
$
x\in G^{\;\!\prime}
\hfill
$
\\[1.75ex]
is orthogonal on the domain $G^{\;\!\prime}$ to the vector fields
\\[1.75ex]
\mbox{}\hfill
$
\bigl(0,\ldots ,0,-\;\!f_{h}(x), 0,\ldots ,0, f_{q}(x),0,\ldots,0\bigr)$
\ for all $x\in G^{\;\!\prime},
\quad 1\leq q < h\leq n,
\hfill
$
\\[1.75ex]
which are associated with differential forms $\psi_{{}_{\scriptstyle qh}}.$
\vspace{0.5ex}

Therefore we shall say that Theotem 2.9 is {\sl the second boundedness test for the number of compact integral hy\-per\-sur\-fa\-ces
of an ordinary autonomous differential system.}
\vspace{0.5ex}

If $n = 2,$ then Theorem 2.9 correspond to Theorem 0.1.
In this case is enough $P\in C^1(G)$ and $Q\in C^1(G),$
since limit cycles, which are closed (rather than only compact) trajectories, are determined by periodic solutions of system (0.1).

Suppose $n = 2,$ the domain $G$ has the fundamental group $\pi_1^{}(G)$ of ranks $0$ and $1.$
Then Theorem 2.9 correspond to the Bendixson --- Dulac test for the absence (when $d(\pi_1^{}(G)) = 0)$
of closed curve, which is made from trajectories of system {\rm (0.1)}, and
correspond to the Dulac test for the existence of at most one (when $d(\pi_{1}^{}(G)) = 1)$
closed curve, which is made from trajectories of system {\rm (0.1)}.
\vspace{0.5ex}

{\bf Theorem 2.10.}
{\it
Let a domain
\vspace{0.35ex}
$G^{\;\!\prime}\subset G$ be an $(r+1)\!$-connected domain, let there exist
scalar functions $\beta\in C^{\;\!1}(G^{\;\!\prime})$ and $\alpha\in C^{\;\!2}(G^{\;\!\prime})$ such that the function
\\[2ex]
\mbox{}\hfill
$
p_{{}_{\scriptstyle 1}}\colon (x,y)\to \
\dfrac{P(x,y)}{Q(x,y)}\ \partial_{x} \alpha(x,y)
$
\ for all $(x,y)\in G^{\;\!\prime}
\hfill
$
\\[2.75ex]
\mbox{}\hfill
$
\biggl( p_{{}_{\scriptstyle 2}}\colon (x,y)\to
\dfrac{Q(x,y)}{P(x,y)}\ \partial_{y} \alpha(x,y)
$
\ for all
$
(x,y)\in G^{\;\!\prime}\biggr)
\hfill
$
\\[2ex]
is continuously differentiable on the domain $G^{\;\!\prime}$ and the function
\\[2ex]
\mbox{}\hfill   
$
q_{{}_{\scriptstyle 1}}\colon (x,y)\to
\partial_{x} \bigl(
p_{{}_{\scriptstyle 1}}(x,y) + \partial_{y} \alpha(x,y)\bigr)
+ {\rm div}\;\!A(x,y)
$
\ for all
$
(x,y)\in G^{\;\!\prime}
$
\hfill {\rm (2.14)}
\\[2.5ex]
\mbox{}\hfill   
$
\Bigl(
q_{{}_{\scriptstyle 2}}\!\colon \!(x,y)\!\to
{}-\partial_y\bigl(\partial_x \alpha(x,y)+
p_{{}_{\scriptstyle 2}}(x,y)\bigr) + {\rm div}\;\!A(x,y)
$
\ for all $(x,y)\in G^{\;\!\prime}
\Bigr)
$
\hfill {\rm (2.15)}
\\[2.25ex]
is constant sign, where the vector field
\\[1.5ex]
\mbox{}\hfill
$
A(x,y)=\beta(x,y)\bigl(P(x,y), Q(x,y)\bigr)
$
\ for all
$
(x,y)\in G^{\;\!\prime}.
\hfill
$
\\[1.5ex]
Then the system {\rm (0.1)} with $P\in C^{\infty}(G)$ and $Q\in C^{\infty}(G)$
in the domain $G^{\;\!\prime}$ has
at most $r$ simple closed curves, which are made from trajectories of system {\rm (0.1)}.
}
\vspace{0.5ex}

{\sl Proof.} Indeed,
\\[2ex]
\mbox{}\hfill
$
d\alpha(x,y)_{\displaystyle |_{(0.1)}} =
\biggl(\;\!\dfrac{P(x,y)}{Q(x,y)}\ \partial_x \alpha(x,y) +
\partial_y\alpha(x,y)\biggr)\;\!dy
$
\ for all $(x,y)\in G^{\;\!\prime},
\hfill
$
\\[2ex]
and
\\[1.25ex]
\mbox{}\hfill
$
d\alpha(x,y)_{\displaystyle |_{(0.1)}} =
\biggl(\;\!\partial_x \alpha(x, y) +
\dfrac{Q(x,y)}{P(x,y)}\ \partial_y \alpha(x,y)\biggr)\;\!dx
$
\ for all
$
(x,y)\in G^{\;\!\prime}.
\hfill
$
\\[2.5ex]
\indent
Then, respectively, the exterior differential
\\[2.5ex]
\mbox{}\hfill
$
d\Bigl(d\alpha(x,y)_{\displaystyle |_{(0.1)}} +
\beta(x,y)\bigl(P(x,y)\,dy - Q(x,y)\,dx\bigr)\Bigr) =
\hfill
$
\\[3ex]
\mbox{}\hfill
$
= \biggl(
\partial_x^{}\biggl(
\dfrac{P(x,y)}{Q(x,y)}\ \partial_x \alpha(x,y) +
\partial_y\alpha(x,y)\biggr) +
\partial_x\bigl(\beta(x,y)P(x,y)\bigr)\, +
\partial_y\bigl(\beta(x,y)Q(x,y)\bigr)\biggr)\;\!dx \wedge dy
\hfill
$
\\[3.25ex]
and
\\[1.75ex]
\mbox{}\hfill
$
d\Bigl(d\alpha(x,y)_{\displaystyle |_{(0.1)}} +
\beta(x,y)\bigl(P(x,y)\,dy - Q(x,y)\,dx\bigr)\Bigr) =
\hfill
$
\\[3ex]
\mbox{}\hfill
$
= \Bigl( {}-\partial_y\Bigl(\partial_x \alpha(x,y) +
\dfrac{Q(x,y)}{P(x,y)}\ \partial_y \alpha(x,y)\Bigr) +
\partial_x\bigl(\beta(x,y)P(x, y)\bigr)  +
\partial_y
\bigl(\beta(x,y)Q(x,y)\bigr)\Bigr)\;\!dx \wedge dy
\hfill
$
\\[2ex]
\mbox{}\hfill
for all
$(x,y) \in G^{\;\!\prime}.
\hfill
$
\\[2.25ex]
\indent
By Theorem 2.3{\rm D}, we conclude that Theorem 2.10 is valid.  \k
\vspace{1ex}

In the case of closed curves, which are limit cycles of system (0.1),
Theorem 2.10 is true for system (0.1) with $P\in C^{\infty}(G)$ and $Q\in C^{\infty}(G),$
since limit cycles, which are closed (rather than only compact) trajectories,
\vspace{0.35ex}
are determined by periodic solutions of system (0.1).

Notice that if under the assumptions of Theorem 2.10 we set
\\[1.5ex]
\mbox{}\hfill
$
\partial_{x}^{}\;\! \alpha(x,y)= \mu(x,y)\;\!Q(x,y)$ \ for all $(x,y) \in G^{\;\!\prime}
\hfill
$
\\[1.5ex]
in the expression (2.14) and
\\[1.5ex]
\mbox{}\hfill
$
\partial_{y}^{}\;\! \alpha(x,y)={}- \mu(x,y)\;\!P(x,y)
$
\ for all
$
(x,y) \in G^{\;\!\prime}
\hfill
$
\\[1.5ex]
in the expression (2.15), then these expressions acquire the form
\\[1.5ex]
\mbox{}\hfill
$
\partial_{x}^{} \bigl( (\mu(x,y)+\beta(x,y))\;\!P(x,y)\bigr)+ \partial_{y}^{}
\bigl( ( \mu(x,y)+\beta(x,y))\;\!Q(x,y)\bigr)
$
\ for all
$(x,y) \in G^{\;\!\prime}.
\hfill
$
\\[2.75ex]
\indent
{\bf 2.3.3.\!\! Test of the boundedness of the number of compact regular integral hy\-per\-sur\-fa\-ces.}
\\[0.35ex]
\indent
{\bf Theorem 2.11}
({\sl a boundedness test for the number of isolated compact regular integral hypersurfaces of
an ordinary autonomous differential system}).
{\it
Let a domain $G^{\;\!\prime}\subset G$ have the homotopy group
$\pi_{n-1}^{}(G^{\;\!\prime})$ of rank $r,$
\vspace{0.5ex}
let there exists a definite on the domain $G^{\;\!\prime}$ scalar function $g\in C^{\infty}(G^{\;\!\prime})$
such that the vector field
\\[1.5ex]
\mbox{}\hfill
$
h\colon x\to g(x) f(x)$
\ for all $x \in G^{\;\!\prime}
\hfill
$
\\[1.5ex]
is solenoidal on the domain $G^{\;\!\prime}.$
\vspace{0.25ex}
Then the system {\rm (D)} with  $f\in C^{\infty}(G)$
has at most $r$ isolated compact regular integral hypersurfaces in the domain $G^{\;\!\prime}.$
}
\vspace{0.5ex}

The {\sl proof} of this theorem is consistent with the proof of Theorem 1.1 when $\nu = n$ and is based on the following
\vspace{0.5ex}

{\bf Lemma 2.3.}
{\it
Let the assumptions of Theorem {\rm 2.11} be valid.
Then the situation described in Lemma {\rm 2.1}
\vspace{0.25ex}
is impossible for isolated compact regular integral hypersurfaces of
the ordinary autonomous differential system {\rm (D)} with $f\in C^{\infty}(G)$
\vspace{0.25ex}
in any subdomain $\Omega$ of the domain $G^{\;\!\prime}\subset G$
with homotopy group $\pi_{n-1}^{}(\Omega)$ of rank $s, \ s \leq r.$
}
\vspace{0.5ex}

{\sl Proof}.
First we note the following fact. Since the function $g\in C^{\infty}(G^{\;\!\prime})$ is definite on the domain $G^{\;\!\prime},$
it follows that if system (D) has a compact regular integral hypersurface in the domain $G^{\;\!\prime},$
then the autonomous ordinary differential system determining the vector field
$
h\colon x\to g(x) f(x)$ for all $x \in G^{\;\!\prime}
$
has the same compact regular integral hypersurface.
\vspace{0.35ex}

Suppose the contrary. Namely, let the situation described in Lemma 2.1 take place for isolated compact regular integral hypersurfaces
$\partial\Sigma_{1}^{},\ldots,\partial\Sigma_{s+1}^{}$ of system {\rm (D)}.
\vspace{0.35ex}
Since $\partial \Sigma_{s+1}^{}$ is an isolated compact regular integral hypersurface of system (D), we see that,
\vspace{0.25ex}
in the exterior, trajectories of system (D) approach $\partial\Sigma_{s+1}^{}$
\vspace{0.35ex}
(respectively, go away from $\partial\Sigma_{s+1}^{})$ as $t \to {}+\infty,$
and there exists a hypersurface $\partial \Xi$ diffeomorphic to the hypersurface $\partial \Sigma_{s+1}^{}$
\vspace{0.25ex}
such that through this hypersurface the trajectories enter
\vspace{0.25ex}
(respectively, leave) the domain bounded by the hypersurfaces $\partial \Xi$ and $\partial \Sigma_{s+1}.$
Therefore,
\\[1.75ex]
\mbox{}\hfill
$
\displaystyle
\int\limits_{\partial \Xi} g(x)f(x) n(x)\,dS \ +\
\sum\limits_{\tau = 1}^s\ \
\int\limits_{\partial \Sigma_{{}_{\scriptstyle \tau}}}
g(x)f(x) n(x)\,dS \ne 0,
\hfill
$
\\[1.75ex]
where n is the unit outward normal field.

But, on the other hand, since the vector field
$
h\colon x\to g(x) f(x)$ for all $x \in G^{\;\!\prime}
$
is solenoidal on the domain $G^{\;\!\prime},$
it follows that
\\[2ex]
\mbox{}\hfill
$
\displaystyle
\int\limits_{\partial \Xi} g(x)f(x) n(x)\,dS \ +\
\sum\limits_{\tau = 1}^s\ \
\int\limits_{\partial \Sigma_{{}_{\scriptstyle \tau}}}
g(x)f(x) n(x)\,dS \ =\
\int\limits_\Xi {\rm div}\;\!\bigl(g(x)f(x)\bigr)\;\!d\Xi \, =\,  0,
\hfill
$
\\[1.75ex]
where $\Xi$ is the domain bounded by the hypersurface
$\partial\Xi\cup\Bigl(\,\bigcup\limits^{s}_{\tau = 1}
\partial \Sigma_{{}_{\scriptstyle \tau}}\Bigr).$
\vspace{1ex}

The obtained contradiction completes the proof of Lemma 2.3. \k

\newpage

From Theorems 2.9 and 2.11 (with $\varphi(x) = 1$ for all $x\in G),$ we obtain the following
\vspace{0.5ex}

{\bf Corollary 2.5.}
{\it
An ordinary linear autonomous differential system has no isolated compact regular integral hypersurface.
}
\\[3.75ex]
\centerline{
\bf 2.4. Autonomous total differential system
}
\\[1.5ex]
\indent
The autonomous system of total differential equations {\rm (TD)}
induces $m$ ordinary autonomous differential systems of $n\!$-th order
{\rm (Dj)}, $j=1,\ldots,m.$
\vspace{0.35ex}

A scalar function
\vspace{0.25ex}
$w\in C^{1}(G^{\;\!\prime}),\ G^{\;\!\prime}\subset G,$ is an autonomous partial integral of system {\rm (TD)}
if and only if the system of identities
\\[1.75ex]
\mbox{}\hfill                                       
$
{\frak X}_{{}_{\scriptstyle j}}w(x) = \Phi_{j}^{}(x)
$
\ for all
$
x \in G^{\;\!\prime},
\quad
j= 1,\ldots,m,
$
\hfill (2.16)
\\[2ex]
is consistent, where the functions $\Phi_j^{}\colon G^{\;\!\prime}\to \R, \ j=1,\ldots,m,$ such that
\\[2.75ex]
\mbox{}\hfill                                       
$
\Phi_{j}(x)_{
\displaystyle |_{\scriptstyle w(x)=0}} = 0,
\quad
j=1,\ldots,m.
$
\hfill (2.17)
\\[2.25ex]
\indent
The consistency of the $k\!$-th identity of system (2.16) under the condition (2.17) with $j=k$
is equivalent to the existence of the autonomous partial integral $w\colon G^{\;\!\prime} \to \R$ for
system {\rm (Dk)}. This allows one to make the following conclusions [9].
\vspace{0.75ex}

{\bf Theorem 2.12}
({\sl
the first boundedness test for the number of compact integral hypersurfaces of an autonomous total differential system}).
{\it
\!Suppose there exists an index $j\in\! \{ 1,\ldots , m\}$ such that the assumptions of Theorem {\rm 2.3D} are valid for system {\rm (Dj)}.
Then the system {\rm (TD)} with $X\in C^{\infty}(G)$
has at most $r$ compact integral hypersurfaces in the domain $G^{\;\!\prime}.$
}
\vspace{0.75ex}

{\bf Theorem 2.13.}
{\it
Suppose there exists an index
\vspace{0.25ex}
$j\in\{1,\ldots,m\}$ such that the assumptions of Theorem {\rm 2.4D} are valid for system {\rm (Dj)}.
Then the following assertions hold{\rm:}
\vspace{0.25ex}

{\rm 1)}
\vspace{0.15ex}
the total differential systen {\rm (TD)} with $X\in C^{\infty}(G)$ has at most $r$ compact integral hypersurfaces in the domain $G^{\;\!\prime};$
\vspace{0.25ex}

{\rm 2)}
\vspace{0.25ex}
any set of gaps of the domain $G^{\;\!\prime}$ surrounded by a compact integral hypersurface of
system {\rm (TD)} with $X\in C^{\infty}(G)$
\vspace{0.75ex}
has the zero total index with respect to the $(n-1)\!$-form $d\theta.$
}

{\bf Theorem 2.14}
({\sl the second boundedness test for the number of compact integral hypersurfaces
of an autonomous total differential system}).
{\it
Suppose a domain $G^{\;\!\prime}\subset G$ has the homotopy group $\pi_{n-1}^{}(G^{\;\!\prime})$ of rank $r,$
\vspace{0.35ex}
there exist an index $j \in \{1,\ldots, m\}$ and a scalar function $\varphi\in C^{1}(G^{\;\!\prime})$
such that the divergence ${\rm div}\,Z_{j}^{}$ of the vector field
\\[1.35ex]
\mbox{}\hfill
$
Z_{j}\colon x \to\,  \varphi(x)\;\!X^j(x)
$
\ for all
$
x \in G^{\;\!\prime}
\hfill
$
\\[1.5ex]
is constant sign on the domain $G^{\;\!\prime}.$ Then the total differential system {\rm (TD)} with $X\in C^{\infty}(G)$
has at most $r$ compact integral hypersurfaces in the domain $G^{\;\!\prime}.$
}
\vspace{1.25ex}

{\bf Example 2.7.}
Consider the total differential system
\\[2ex]
\mbox{}\hfill                                         
$
\displaystyle
dx_1^{} = x_3^{}\prod\limits_{k=0}^{n}\,(g(x) - 2k)\,dt_1^{} + x_2^{}\,dt_2^{}\;\!,
\hfill
$
\\[2.25ex]
\mbox{}\hfill
$
\displaystyle
dx_2^{} = \biggl(x_3^{} + x_2^{}\prod\limits_{k=0}^{n}\,
(g(x) - 2k)(g(x) - 2k - 1)\biggr)\;\!dt_1^{} + ({}-x_1^{} + x_3^{})\,dt_2^{}\;\!,
$
\hfill (2.18)
\\[2.25ex]
\mbox{}\hfill
$
\displaystyle
dx_3^{} =
\biggl({}-x_2^{} - x_1^{}\prod\limits_{k=0}^{n} (g(x) - 2k)\,
\biggr)\;\!dt_1^{}  +
\biggl({}-x_2^{} + x_3^{}\prod\limits_{k=0}^{n}\,
(g(x) - 2k)(g(x) - 2k - 1)\biggr)\;\!dt_2^{}\;\!,
\hfill
$
\\[2ex]
where the scalar function
\\[0.5ex]
\mbox{}\hfill
$
g\colon x\to\,  x_1^2 + x_2^2 + x_3^2$ \ for all $x\in\R^3.
\hfill
$
\\[1.5ex]
\indent
Since
\\[1.5ex]
\mbox{}\hfill
$
\displaystyle
d\bigl(g(x) - m\bigr)_{\displaystyle |_{(2.18)}} =
2\prod\limits_{k=0}^{n}\,(g(x) - 2k)(g(x) - 2k - 1)(x_2^2\,dt_1^{} + x_3^2\,dt_2^{})
\hfill
$
\\[2ex]
\mbox{}\hfill
for all
$
(t, x) \in \R^5,
\quad
m = 1,\ldots, 2n +1,
\hfill
$
\\[1.75ex]
we see that the spheres
\vspace{0.35ex}
$
S^2_{m} = \{x\colon g(x) = m\}, \
m = 1,\ldots, 2n +1,
$
are compact integral hy\-per\-sur\-fa\-ces for the system of total differential equations (2.18).
\vspace{0.15ex}

Using the ordinary autonomous differential system {\rm (D1)} induced by the
autonomous total differential system (2.18), we obtain the Pfaffian equations
\\[1.75ex]
\mbox{}\hfill
$
\psi_{{}_{\scriptstyle 12}}(x) =0,
\hfill
$
\\[1.25ex]
where the differential 1-form
\\[2.25ex]
\mbox{}\hfill
$
\displaystyle
\psi_{{}_{\scriptstyle 12}}(x) =
{}-\biggl(x_3^{} + x_2^{}\prod\limits_{k=0}^{n}\,(g(x) - 2k)(g(x) - 2k - 1)\biggr)\;\!dx_1^{} \,+
x_3^{}\prod\limits_{k=0}^{n} \,(g(x) - 2k)\,dx_2^{}
\hfill
$
\\[1.75ex]
\mbox{}\hfill
for all $x \in\R^3.
\hfill
$
\\[1.25ex]
\indent
The linear differential form
\\[2.25ex]
\mbox{}\hfill
$
\displaystyle
\ell(x) ={}- \prod\limits_{k=0}^{n} \
\dfrac{1}{g(x)-2k}\  dx_1^{}
$
\ for all $x \in \bigcup\limits^{n}_{k=0} G_k,
\hfill
$
\\[2.5ex]
\mbox{}\hfill
$
\displaystyle
G_{{}_{\scriptstyle \tau}} =
\{x\colon 2\tau < g(x) < 2(\tau +1)\},
\ \tau= 0,\ldots,n-1,
\quad
G_{{}_{\scriptstyle n}} = \{x\colon x_1^2 + x_2^2 + x_3^2 > 2n\},
\hfill
$
\\[2.25ex]
such that the exterior differential
\\[2.25ex]
\mbox{}\hfill
$
\displaystyle
d\bigl(\psi_{{}_{\scriptstyle 12}}(x) \wedge \ell(x)\bigr) =
dx_1\wedge
dx_2\wedge
dx_3
$
\ for all
$x \in \bigcup\limits^{n}_{k=0} G_{{}_{\scriptstyle k}}.
\hfill
$
\\[2ex]
\indent
Therefore, by Theorem 2.12, the system (2.18) in any domain
\vspace{0.25ex}
$G_k^{}$ with homotopy group $\pi_{2}^{}(G_k^{}),\ k=0,\ldots, n,$ has
at most one compact integral hypersurface.
\vspace{0.5ex}

Thus the autonomous total differential system (2.18) has
\vspace{0.15ex}
$2n+1$ compact integral hypersurfaces, which are the spheres $S_m^2,\ m = 1,\ldots, 2n +1.$
\vspace{0.75ex}

If system {\rm(TD)} is completely solvable, then in the entire set of its integral hypersurfaces we
distinguish regular ones, on which this system has no singular points.
\vspace{0.75ex}

{\bf Theorem 2.15}
({\sl a boundedness test for the number of isolated compact regular integral hypersurfaces of
a completely solvable autonomous total differential system}).
\vspace{0.25ex}
{\it
Suppose a domain $G^{\;\!\prime}\subset G$ has the homotopy group $\pi_{n-1}^{}(G^{\;\!\prime})$ of rank $r,$
\vspace{0.35ex}
there exist definite on the domain $G^{\;\!\prime}$
scalar functions
$g_j^{}\in C^{\infty}(G^{\;\!\prime}),\ j = 1,\ldots, m,$
such that the vector fields
\\[1.75ex]
\mbox{}\hfill
$
Y_j^{}\colon x \to\ g_j^{}(x)\;\!X^j(x)
$
\ for all
$
x\in G^{\;\!\prime},
\quad
j = 1,\ldots, m,
\hfill
$
\\[2ex]
are solenoidal on the domain $G^{\;\!\prime}.$
\vspace{0.25ex}
Then the completely solvable system {\rm (TD)} with $X\in C^{\infty}(G)$
has at most $r$ isolated compact regular integral hypersurfaces in the domain $G^{\;\!\prime}.$
}
\vspace{0.5ex}

{\sl Proof.}
The assertion of this theorem follows from the obvious fact that if each ordinary differential system
{\rm (Dj)}, $j=1,\ldots,m,$ has the same isolated compact integral hypersurface, then this hypersurface
is an isolated compact integral hypersurface of system {\rm (TD)}. \k
\vspace{0.75ex}

Under the conditions (2.17) we can obtain the following statement from the system of identities (2.16),
Corollary 2.5, Theorems 2.12  and  2.15.
\vspace{0.5ex}

{\bf Corollary 2.6.}
{\it
A linear completely solvable autonomous total differential system has no isolated compact regular integral hypersurfaces.
}

\newpage

\mbox{}
\\[-1.75ex]
\centerline{
\bf 2.5. Linear homogeneous system of partial differential equations
}
\\[1.5ex]
\indent
The linear homogeneous system of partial differential equations $(\partial\;\!)$
induces $m$ ordinary autonomous differential systems of $n\!$-th order
{\rm (Dj)}, $j=1,\ldots,m$ (the characteristic system for the partial system $(\partial\;\!)).$
Therefore a scalar function
\vspace{0.25ex}
$w\in C^{1}(G^{\;\!\prime}),\ G^{\;\!\prime}\subset G,$ is an partial integral of system $(\partial\;\!)$
if and only if the system of identities (2.16) under the conditions (2.17) is consistent.
The consistency of the $k\!$-th identity of system (2.16) under the condition (2.17) with $j=k$
is equivalent to the existence of the autonomous partial integral $w\colon G^{\;\!\prime} \to \R$ for
system {\rm (Dk)}. This allows one to make the following conclusions.
\vspace{0.75ex}

{\bf Theorem 2.16}
({\sl
the first boundedness test for the number of compact integral hypersurfaces of a
linear homogeneous system of partial differential equations}).\!
{\it
Suppose there exists an index $\!j\!\in\! \{ 1,\ldots , m\}\!$ such that the assumptions of Theorem {\rm 2.3D} are valid for system\! {\rm (Dj)}.
Then the partial system $\!(\partial\;\!)\!$ with ${\frak X}_{{}_{\scriptstyle j}}\in C^{\infty}(G),\, j\!=\!1,\ldots,m,$
has at most $r$ compact in\-teg\-ral hypersurfaces in the domain $G^{\;\!\prime}.$
}
\vspace{0.75ex}

{\bf Theorem 2.17.}
{\it
Suppose there exists an index
\vspace{0.25ex}
$j\in\{1,\ldots,m\}$ such that the assumptions of Theorem {\rm 2.4D} are valid for system {\rm (Dj)}.
Then the following assertions hold{\rm:}
\vspace{0.25ex}

{\rm 1)}
\vspace{0.15ex}
the linear homogeneous system of partial differential equations
$(\partial\;\!)$ with  ${\frak X}_{{}_{\scriptstyle j}}\in C^{\infty}(G),$ $j=1,\ldots,m,$
has at most $r$ compact integral hypersurfaces in the domain $G^{\;\!\prime};$
\vspace{0.35ex}

{\rm 2)}
\vspace{0.25ex}
any set of gaps of the domain $G^{\;\!\prime}$ surrounded by a compact integral hypersurface of
system $(\partial\;\!)$
\vspace{0.75ex}
has the zero total index with respect to the $(n-1)\!$-form $d\theta.$
}

{\bf Theorem 2.18}
({\sl the second boundedness test for the number of compact integral hypersurfaces
of a linear homogeneous system of partial differential equations}).
\vspace{0.25ex}
{\it
Suppose a domain $G^{\;\!\prime}\subset G$ has the homotopy group $\pi_{n-1}^{}(G^{\;\!\prime})$ of rank $r,$
\vspace{0.5ex}
there exist an index $j \in \{1,\ldots, m\}$ and a scalar function $\varphi\in C^{1}(G^{\;\!\prime})$
such that the divergence ${\rm div}\,Z_{j}^{}$ of the vector field
\\[1.5ex]
\mbox{}\hfill
$
Z_{j}\colon x \to\,  \varphi(x)\;\!X^j(x)
$
\ for all
$
x \in G^{\;\!\prime}
\hfill
$
\\[1.75ex]
is constant sign on the domain $G^{\;\!\prime}.$ Then
\vspace{0.25ex}
the system of partial differential equations
$(\partial\;\!)$ with  ${\frak X}_{{}_{\scriptstyle j}}\in C^{\infty}(G),\ j=1,\ldots,m,$
\vspace{1.25ex}
has at most $r$ compact integral hypersurfaces in the domain $G^{\;\!\prime}.$
}

{\bf Example 2.8.}
Consider the linear homogeneous system of partial differential equations
\\[2.25ex]
\mbox{}\hfill
$
x_2^{}g(x)\;\!\partial_{{}_{\scriptstyle x_1^{}}}y +
\bigl(x_3^{} - x_1^{}\;\!g(x)\bigr)\;\!\partial_{{}_{\scriptstyle x_2^{}}}y +
\bigl({}-x_1^{} - x_2^{} + x_1^{}\;\!g(x)\bigr)^{}
\partial_{{}_{\scriptstyle x_3^{}}}y = 0,
\hfill
$
\\
\mbox{}\hfill(2.19)
\\[0.5ex]
\mbox{}\hfill
$
\bigl({}-x_2^{} + x_3^{} + x_2^{}\;\!g(x)\bigr)\;\!
\partial_{{}_{\scriptstyle x_1^{}}} y +
\bigl(x_1^{} - x_1^{}\;\!g(x)\bigr)
\partial_{{}_{\scriptstyle x_2^{}}} y  +
\bigl({}-x_1^{} - x_2^{} + x_2^{}\;\!g(x)\bigr)\;\!
\partial_{{}_{\scriptstyle x_3^{}}} y =  0,
\hfill
$
\\[2ex]
where the scalar function
$
g\colon x\to x_1^2 + x_2^2 + x_3^2$ for all $x\in\R^3.$
\vspace{0.5ex}

The sphere
$
S^2 = \{x\colon g(x) = 1\}
$
is a compact integral surface of system (2.19).
\vspace{0.35ex}

Now let us consider the ordinary autonomous differential system {\rm (D1)} induced by the partial system (2.19).
Using the system {\rm (D1)}, we obtain the Pfaffian equations
\\[1.5ex]
\mbox{}\hfill
$
\psi_{{}_{\scriptstyle 13}}(x) = 0,
\hfill
$
\\[1ex]
where the differential 1-form
\\[1.75ex]
\mbox{}\hfill
$
\psi_{{}_{\scriptstyle 13}}(x) =x_2^{}\;\!g(x)\,dx_3^{} - \bigl({}-x_1^{} - x_2^{} + x_1^{}\;\!g(x)\bigr)\;\!dx_1^{}
$
\ for all
$
x \in \R^3.
\hfill
$
\\[2ex]
\indent
The continuously differentiable on the domain $G^{\;\!\prime} = \R^3\backslash \{(0,0,0)\}$
differential 1-form
\\[2ex]
\mbox{}\hfill
$
\ell(x)=\dfrac{1}{g(x)}\ dx_1^{}
$
\ for all
$
x\in G^{\;\!\prime}
\hfill
$
\\[2ex]
such that the exterior differential of exterior product
\\[2ex]
\mbox{}\hfill
$
d\bigl(\psi_{{}_{\scriptstyle 13}}(x) \wedge \ell(x)\bigr) =
dx_{1}^{} \wedge  dx_{2}^{} \wedge dx_{3}^{}
$
\ for all
$
x \in G^{\;\!\prime}.
\hfill
$
\\[2ex]
\indent
Therefore, by Theorem 2.16, the system (2.19) in the domain
\vspace{0.25ex}
$G^{\;\!\prime}$ with homotopy group $\pi_{2}(G^{\;\!\prime})$ of rank 1
has at most one compact integral surface.
This surface is the sphere $S^2.$
\vspace{1ex}

{\bf Example 2.9.}
For the linear homogeneous system of partial differential equations
\\[2ex]
\mbox{}\ \
$
\bigl({}-x_1^{}+ x_2^{} + x_1^{}\;\!g(x)\bigr)
\partial_{{}_{\scriptstyle x_1^{}}} y +
\bigl({}-x_1^{} - x_2^{} + x_2^{}\;\!g(x)\bigr)
\partial_{{}_{\scriptstyle x_2^{}}}y +
\bigl({}-x_3^{}+ x_3^{}\;\!g(x)\bigr)
\partial_{{}_{\scriptstyle x_3^{}}}y=0,
\hfill
$
\\[0.5ex]
\mbox{}\hfill (2.20)
\\[0.5ex]
\mbox{}\quad \,
$
\bigl({}-x_2^{} + x_3^{} + x_2^{}\;\!g(x)\bigr)
\partial_{{}_{\scriptstyle x_1^{}}}y +
\bigl(x_1^{} - x_1^{}\;\!g(x)\bigr)
\partial_{{}_{\scriptstyle x_2^{}}}y  +
\bigl({}-x_1^{} - x_2^{} + x_2^{}\;\!g(x)\bigr)
\partial_{{}_{\scriptstyle x_3^{}}}y = 0,
\hfill
$
\\[2.25ex]
where the scalar function
$
g\colon x\to x_1^2 + x_2^2 + x_3^2$ for all $x\in\R^3,
$
the sphere
\\[2ex]
\mbox{}\hfill
$
S^2 = \{x\colon g(x) = 1\}
\hfill
$
\\[1ex]
is a compact integral surface.

Consider now the ordinary autonomous differential system {\rm (D1)} induced by the
partial differential system (2.20).
Let us the scalar function
\\[1.5ex]
\mbox{}\hfill
$
\varphi\colon x \to \
g^{^{{}-\tfrac52}}(x)
$
\ for all
$
x\in\R^3\backslash \{(0,0,0)\}.
\hfill
$
\\[2ex]
Then the divergence ${\rm div}\,Z_{1}^{}$ of the vector field
\\[2ex]
\mbox{}\hfill
$
Z_1^{}\colon x \to\ \varphi(x)\;\!X^1(x)
$
\ for all $x\in \R^3\backslash \{(0,0,0)\}
\hfill
$
\\[2ex]
is positive on the domain $\R^3\backslash\{(0,0,0)\}.$
The rank
$
d(\pi_{2}^{}(\R^3\backslash\{(0,0,0)\}))=1.
$
\vspace{0.5ex}

By Theorem 2.18, it follows that
\vspace{0.35ex}
the sphere $S^2$ is the unique compact integral surface of
the system of partial differential equations (2.20)
in the domain $\R^3\backslash \{(0,0,0)\}.$

}
\end{document}